\documentclass[11pt]{amsart}
\usepackage{amsfonts,latexsym,amsthm,amssymb,amsmath,amscd,euscript,tikz, tikz-cd}
\usepackage[alphabetic, msc-links, bibtex-style, nobysame]{amsrefs}
\usepackage{stackengine}
\usepackage{framed}
\usepackage{xfrac}
\usepackage[makeroom]{cancel}
\usepackage{faktor}
\usepackage{braket}
\usepackage{pgf,tikz,pgfplots}
\pgfplotsset{compat=1.17}
\usepgfplotslibrary{fillbetween} 
\usepackage{mathrsfs}
\usetikzlibrary{arrows}
\definecolor{wrwrwr}{rgb}{0.3803921568627451,0.3803921568627451,0.3803921568627451}
\definecolor{rvwvcq}{rgb}{0.08235294117647059,0.396078431372549,0.7529411764705882}
\definecolor{mblue}{rgb}{0.2, 0.3, 0.8}
\definecolor{morange}{rgb}{1, 0.5, 0}
\definecolor{mgreen}{rgb}{0.1, 0.4, 0.2}
\definecolor{mred}{rgb}{0.5, 0, 0}
\definecolor{ForestGreen}{RGB}{34,139,34}
\usepackage{hyperref}
\hypersetup{colorlinks=true,citecolor=ForestGreen,linkcolor = blue,urlcolor =black,linkbordercolor={1 0 0}}
\usepackage{float}

\usepackage{lmodern}
\usepackage{supertabular}
\usepackage{verbatim}
\usepackage{amssymb}
\usepackage{enumerate}
\usepackage{stmaryrd}
\usepackage{bbm}
\usepackage{mathtools}
\usepackage[nopatch=footnote]{microtype}
\usepackage{setspace}
\usepackage{tikz,tikz-cd}
\usetikzlibrary{matrix,calc,positioning,arrows,decorations.pathreplacing,patterns,knots}

\newcommand{\la}{\langle}
\newcommand{\rg}{\rangle}

\newtheorem{theorem}{{Theorem}}[section]
\newtheorem*{theorem*}{Theorem}
\newtheorem{lemma}[theorem]{Lemma}
\newtheorem{proposition}[theorem]{Proposition}

\newtheorem{corollary}[theorem]{Corollary}
\newtheorem*{corollary*}{Corollary}

\theoremstyle{definition}

\usepackage{lmodern,url,enumerate,mathtools,microtype}
\usepackage[hmargin = 1in,vmargin=1in]{geometry}
\usepackage{graphicx}
\usepackage{subcaption}

\newcommand{\ve}{\varepsilon}

\newcommand{\mr}[1]{{\rm #1}}
\newcommand{\cC}{\mathcal{C}}

\newcommand{\cH}{\mathcal{H}}

\newcommand{\cL}{\mathcal{L}}

\newcommand{\cR}{\mathcal{R}}

\newcommand{\bR}{\mathbb{R}}
\newcommand{\bS}{\mathbb{S}}

\newcommand{\nc}{\newcommand}

\nc{\on}{\operatorname}
\nc{\p}{\partial}
\nc{\ol}{\overline}
\nc{\ul}{\underline}
\nc{\pa}{\partial}

\nc{\pb}{\partial_b}
\nc{\pc}{\partial_c}
\nc{\pd}{\partial_d}
\nc{\pe}{\partial_e}
\nc{\pf}{\partial_f}
\nc{\pg}{\partial_g}
\nc{\ph}{\partial_h}
\nc{\pari}{\partial_i}
\nc{\pj}{\partial_j}
\nc{\pk}{\partial_k}
\nc{\pl}{\partial_l}
\nc{\pell}{\partial_\ell}
\nc{\parm}{\partial_m}
\nc{\pn}{\partial_n}
\nc{\po}{\partial_o}
\nc{\pp}{\partial_p}
\nc{\pq}{\partial_q}
\nc{\pr}{\partial_r}
\nc{\ps}{\partial_s}
\nc{\pt}{\partial_t}
\nc{\pu}{\partial_u}
\nc{\pv}{\partial_v}
\nc{\pw}{\partial_w}
\nc{\px}{\partial_x}
\nc{\py}{\partial_x}
\nc{\pz}{\partial_z}

\nc{\Spec}{\on{Spec}}

\nc{\sn}{\mr{sn}}
\nc{\cn}{\mr{cn}}
\nc{\dn}{\mr{dn}}

\DeclareMathOperator*{\essinf}{ess\,inf}
\DeclareMathOperator*{\esssup}{ess\,sup}

\allowdisplaybreaks

\title[Stable $s$-minimal surfaces]{Stable $s$-minimal hypersurfaces \\ in dimensions seven and eight}
\date{\today}

\author[Serena Dipierro]{Serena Dipierro}
\address{
{\vspace{-0.1in} \newline
Department of Mathematics and Statistics, The University of Western Australia
\newline \href{mailto:serena.dipierro@uwa.edu.au}{serena.dipierro@uwa.edu.au}} }
\author[Joaquim Serra]{Joaquim Serra}
\address{ \vspace{-0.1in} \newline 
Department of Mathematics, ETH Zurich
\newline {\href{mailto:joaquim.serra@math.ethz.ch}{joaquim.serra@math.ethz.ch}}}
\author[Raphael Tsiamis]{Raphael Tsiamis}
\address{ \vspace{-0.1in} \newline 
Department of Mathematics, Columbia University
\newline {\href{mailto:r.tsiamis@columbia.edu}{r.tsiamis@columbia.edu}}}
\author[Enrico Valdinoci]{Enrico Valdinoci}
\address{
{ \vspace{-0.1in} \newline Department of Mathematics and Statistics, The University of Western Australia \newline \href{mailto:enrico.valdinoci@uwa.edu.au}{enrico.valdinoci@uwa.edu.au}}}
\begin{document}

\begin{abstract}
We establish optimal sheet separation estimates for stable $s$-minimal hypersurfaces in dimensions up to seven, using an iteration argument in the nonlocal approximation theory of Chan-Dipierro-Serra-Valdinoci. These estimates imply that every stable $s$-minimal hypersurface in $\mathbb{R}^7$ is a hyperplane when $s$ is close to $1$. More generally, every finite Morse index $s$-minimal set in $\mathbb{R}^7$ is a half-space.

We also establish a compactness theorem for stable $s$-minimal cones in $\mathbb{R}^8$. We prove uniform curvature and area estimates for their links and show that, as $s \uparrow 1$, they subconverge graphically in $C^{2,\alpha}$, either with multiplicity one to a classical stable minimal cone, or with multiplicity two to a hyperplane. The stability of Lawson and axisymmetric $s$-minimal cones in $\mathbb{R}^8$ shows that these results are optimal.
\end{abstract}

\vspace*{-0.2in}
\maketitle

\vspace*{-0.1in}

\section{Introduction}

Nonlocal minimal surfaces arise as minimizers of the fractional perimeter, a geometric energy that accounts for long-range interactions across an interface. 
Their theory was initiated by Caffarelli, Roquejoffre, and Savin~\cite{caffarelli-roquejoffre-savin}, providing a natural nonlocal analogue of classical minimal surfaces that appears in models involving phase transitions, anomalous diffusion, and interacting particle systems.
In the near-local regime of $s$ sufficiently close to $1$, Caffarelli-Valdinoci~\cite{caffarelli-valdinoci-advances} established the flatness of minimizing cones in $\bR^n$, for every $2 \leq n \leq 7$.
The study of stable (instead\footnote{As is customary in the theory of nonlocal minimal surfaces, we reserve the term {\em$s$-minimal sets} (or {\em minimizing minimal surfaces}) for local minimizers of the~\(s\)-perimeter functional. Correspondingly, we refer to critical points of the~\(s\)-perimeter with nonnegative second variation with respect to smooth vector fields as {\em stable minimal surfaces}. See~\cite{nonlocal-approximation}*{pages~3--5} for detailed notations and definitions.} of minimizing)
nonlocal minimal surfaces in the near-local regime was developed significantly by Cabr\'e-Cinti-Serra~\cites{cabre-cinti-serra , s-allen-cahn}, who proved the flatness of stable $s$-minimal cones in $\mathbb{R}^3$ and used this result to obtain fractional versions of the De Giorgi conjecture and the stable Bernstein theorem.
Previously, Farina, Figalli, and Valdinoci~\cites{flatness-farina-valdinoci , figalli-valdinoci } showed that the flatness of $s$-minimizing cones in $\bR^n$ implies the stable nonlocal Bernstein theorem in $\mathbb{R}^n$, and the graphical Bernstein theorem in $\bR^{n+1}$~\cites{figalli-valdinoci , flatness-farina-valdinoci}.

Chan-Dipierro-Serra-Valdinoci~\cite{nonlocal-approximation} developed a robust approximation theory as $s\uparrow1$ and proved optimal $L^{-2}$ separation estimates between collapsing sheets, namely an $L^2$ estimate on the reciprocal of the distance between sheets. 
Additional results in~\cite{nonlocal-approximation} establish
a suitable decoupling of the $s$-minimal surface equation along different sheets and sharp regularity estimates once pointwise separation is available.
These ingredients lead to curvature estimates using a blow-up argument, in which the limiting object is a stable classical minimal hypersurface.
Using this new theory, \cite{nonlocal-approximation} showed that stable $s$-minimal hypersurfaces in $\mathbb{R}^4$ are flat for $s$ close to $1$.

Notably, the results in~\cite{nonlocal-approximation} did not require the corresponding stable Bernstein theorem in $\mathbb{R}^4$, which was proved by Chodosh and Li~\cite{chodosh-li-R4}, and, later and using alternative methods, by~\cites{chodosh-li , catino-mastrolia-roncoroni , hong-wang}.
In their work, Catino, Mastrolia, and Roncoroni~\cite{catino-mastrolia-roncoroni} also showed that there are no complete orientable stable minimal immersions $M^n \looparrowright \bS^{n+1}$ for $n \leq 5$, a result recently extended to minimal immersions into $\bS^7$ by~\cite{wang-bernstein}.

Since the paper~\cite{nonlocal-approximation} originally appeared, the stable Bernstein problem in $\mathbb{R}^n$ has been completely resolved, thanks to the groundbreaking works of Chodosh-Li-Minter-Stryker in $\mathbb{R}^5$~\cite{chodosh-li-minter-stryker}, Mazet in $\mathbb{R}^6$~\cite{mazet}, and, most recently, by Hong-Li-Wang in $\mathbb{R}^7$~\cite{hong-li-wang}.
The latter breakthrough used the Green's function approach described in~\cite{cabre-catino-gauss}.
This result also relied on the important work of Bellettini, which refined~\cites{schoen-simon-yau , sheeting-theorem } using a De Giorgi-type iteration~\cite{bellettini}.

We show that the approximation theory and the estimates developed in~\cite{nonlocal-approximation} can be iterated in the spirit of De Giorgi and Schoen-Simon-Yau.
Combining this result with the recent advances on stable minimal surfaces~\cites{chodosh-li-minter-stryker , mazet , hong-li-wang  }, we prove the flatness of stable $s$-minimal hypersurfaces in ambient dimension at most $7$. In this framework, our main result reads as follows:

\begin{theorem}\label{thm:stable-cones-s-1}
There exists $s_* \in(0, 1)$ such that, for every $s \in (s_*, 1)$ and every $3 \leq n \leq 7$, every stable $s$-minimal cone in $\bR^n$ with $C^2$ boundary away from the origin is a half-space.
\end{theorem}
In $\mathbb{R}^8$, the Simons cone $\{ (x,y) \in \bR^4 \times \bR^4 : |x| >|y| \} $ is $s$-minimal and stable for $s$ close to $1$, by the computation in~\cite{davila-delpino-wei}.
Therefore, Theorem~\ref{thm:stable-cones-s-1} determines the sharp regularity dimension $n^*_s$ for which there are non-flat stable $s$-minimal cones as $n_s^* = 8$, when $s \in (s_*, 1)$.

In eight dimensions, we obtain a structural theorem for stable $s$-minimal cones, which entails uniform curvature and area estimates as $s \uparrow 1$.
In what follows, we denote by $A_{\Sigma}$ the second fundamental form of a hypersurface $\Sigma$, and write $|A_{\Sigma}|$ for its pointwise norm with respect to the induced metric on $\Sigma$.
\begin{theorem}\label{thm:eight-dimensions}
There exist $s_* \in(0, 1)$ and $C\in(0,\infty)$ such that, for $s \in (s_*,1)$, every stable $s$-minimal cone $E \subset \mathbb{R}^8$ with vertex at the origin and $C^2$ link satisfies
\[
\sup_{\partial E \cap \bS^7} |A_{\partial E \cap \bS^7}| + \cH^6(\partial E \cap \bS^7) \leq C.
\]
Moreover, given a sequence of $s_k \in (s_*,1)$ and stable $s_k$-minimal cones $E_k$, after passing to a subsequence $s_k \to s_{\infty} \in [s_0,1]$, one of the following holds.
\begin{enumerate}[(i)]
    \item If $s_{\infty}<1$, the links converge with multiplicity one to the $C^2$ link of a stable $s_{\infty}$-minimal cone.
    \item If $s_{\infty}=1$, the links converge either with multiplicity one to the link of a non-flat stable classical minimal cone, or with multiplicity two to an equator.
    For a sequence of non-flat cones, the multiplicity-one limit is non-flat.
\end{enumerate}
In all cases, the convergence is graphical in $C^{2,\gamma}$, for every given $\gamma \in (0, \frac{1}{16})$.    
\end{theorem}
In particular, Theorem~\ref{thm:eight-dimensions} yields a compactness theorem for stable $s$-minimal cones in $\mathbb{R}^8$ as $s \uparrow 1$, which is novel even in the case of classical minimal cones (see the forthcoming
Theorem~\ref{thm:density-classical-cones} for the statement in the classical case).
We note that multiplicity-two convergence cannot be ruled out in case $(ii)$, and indeed occurs for the axisymmetric cones considered in~\cite{davila-delpino-wei}.
These cones have the form
\[
E_s = \{ (x', x_8) \in \bR^7 \times \bR : |x_8| = \alpha_s |x'| \}, \qquad \alpha_s = \tfrac{1}{\sqrt{6}} \sqrt{1-s} + o ( \sqrt{1-s}) \quad \text{as } \; s \uparrow 1,
\]
are stable, and converge with multiplicity two to the hyperplane $\bR^7 \times \{ 0 \}$.

Thanks to the works~\cites{figalli-valdinoci , flatness-farina-valdinoci , cabre-cinti-serra , s-allen-cahn }, Theorem~\ref{thm:stable-cones-s-1} has further important consequences to the Bernstein and De Giorgi problems.
Caselli, Florit-Simon, and Serra~\cite{caselli-fs-serra} generalized the above flatness properties to solutions with finite Morse index and proved a nonlocal version of Yau's conjecture for $s$-minimal surfaces.
Combining these results, we obtain:
\begin{corollary}\label{cor:all-the-corollaries}
There exists $s_* \in (0,1)$ such that, for every~$s \in (s_*, 1)$ and $3 \leq n \leq 7$, the following properties hold.
\begin{enumerate}[(i)]
    \item Every $s$-fractional Allen-Cahn solution in $\bR^n$ with finite Morse index is one-dimensional (meaning that it depends on a single Cartesian variable up to a spatial rotation).
    \item Every entire $C^2$ $s$-minimal surface in $\bR^n$ with finite Morse index is a half-space.
    \item Every closed $n$-manifold contains infinitely many smooth $s$-minimal surfaces.
    More generally, in dimensions $n \geq 8$, the singular set has Hausdorff dimension less than or equal to~$n-8$, with a discrete singular set in dimension $8$.
    \item Every entire globally minimizing $s$-minimal graph in $\bR^8$ is a hyperplane.
\end{enumerate}
\end{corollary}
We note that the strongly nonlocal regime $s \sim 0$ is less well understood.
See~\cite{tsiamis} for the resolution of the corresponding problems in $\mathbb{R}^3$.

We obtain Theorems~\ref{thm:stable-cones-s-1} and~\ref{thm:sheet-happens} through an improvement of the $L^{-2}$ sheet separation estimate, which was central in~\cite{nonlocal-approximation}.
Notably, in the forthcoming Proposition~\ref{prop:improve-the-exponent}, we establish an $L^{-6}$ sheet separation estimate for stable $s$-minimal hypersurfaces with $s \sim 1$, which implies the following structural result.
\begin{theorem}\label{thm:sheet-happens}
There exist constants $\delta_* > 0$, $s_* \in (0,1)$, $r_* > 0$, $c>0$, and $C \in(0, \infty)$, with the following property.

Let $2 \leq n \leq 7$ and $s \in (s_*, 1)$. Let $E \subset \bR^n$ be a stable $s$-minimal set in the cylinder $\Omega := B_{1/2}^{n-1} \times (-1,1) \subset \bR^n$.
Suppose that the boundary $\partial E$ in $\Omega$ is a union of graphs 
\[
\partial E \cap \Omega = \bigcup_{j=1}^N \{ (x', g_j(x')) : x' \in B_{1/2}^{n-1}, \; -1<g_j(x') < 1 \}
\]
inside $\Omega$, for functions $g_j \in C^2(B^{n-1}_1)$ that satisfy
\[
g_1<g_2<\cdots< g_N \qquad{\mbox{and}} \qquad\sup_j \, ( \|D g_j \|_{L^{\infty}(B^{n-1}_1)} + \|D^2 g_j \|_{L^{\infty}(B_1^{n-1})}) < \delta_*.
\]
Then, every adjacent pair of graphs obeying the bounds $|g_i(0)| \leq \frac{1}{4}$ and $|g_{i+1}(0)| \leq \frac{1}{4}$ satisfies
\[
\inf_{B_{r_*}^{n-1}} (g_{i+1} - g_i) \geq c \, \sqrt{1-s}.
\]
Moreover, the graph $\Gamma_i = \{ x_n = g_i(x') \}$ of every sheet with $|g_i(0)| \leq \frac{1}{4}$ satisfies the uniform estimates
\[
\| g_i \|_{C^{2,1/8}(B^{n-1}_{r_*/4})} \leq C \qquad{\mbox{and}} \qquad \sup_{x' \in B^{n-1}_{r_*/4}} \bigl| H[\Gamma_i](x', g_i(x')) \bigr| \leq C \sqrt{1-s}.
\]
The constants are independent of the number $N$ of sheets.

The same conclusion also holds in $\mathbb{R}^8$ if the stable $s$-minimal set $E$ is, in addition, conical about some point $x_{\circ}$ with $|x_{\circ}| \geq 100$ (meaning that $t (E - x_{\circ}) = E - x_{\circ}$ for every $t>0$).
\end{theorem}
This result extends the sheet separation and regularity theory of~\cite{nonlocal-approximation} to higher dimensions, providing a nonlocal analogue of the corresponding estimates for Allen-Cahn approximations established by Chodosh-Mantoulidis~\cite{chodosh-mantoulidis} and Wang-Wei~\cite{second-order-estimates}.

\subsection{Strategy of the proof}
Our results proceed in four steps, following the theory of~\cite{nonlocal-approximation}.

In Section~\ref{sec:prelimiinaries-nonlocal-approximation}, we establish some important preliminaries from the work~\cite{nonlocal-approximation}, notably the construction of a nonlocal stability operator.
In Lemma~\ref{lemma:slab-estimate}, we control the contribution from exterior terms and interactions of distant sheets to the stability inequality.
Then, in Proposition~\ref{prop:stable-toda}, we establish an $L^{-p}$ estimate, for any $p < 4 + 2 \sqrt{2}$, for stable solutions of the limiting problem.
This result is closely related to the Schoen-Simon-Yau iteration for the second fundamental form~\cite{schoen-simon-yau}, and extends the work of Meadows for the Lane-Emden equation $\Delta u = \frac{1}{u}$, cf.~\cite{meadows}*{\S 7 - 9}.
We refer the reader to the discussion of~\cite{nonlocal-approximation}*{\S 1.8} regarding the relationship between this problem and the system satisfied by the sheets in the near-local regime $s \uparrow 1$.
While this result is not strictly necessary for the ensuing arguments, it elucidates the argument we follow in the nonlocal case.

For our key result, in Proposition~\ref{prop:improve-the-exponent}, we implement a nonlocal iteration to obtain an $L^{-6}$ estimate for the sheet separation of stable $s$-minimal sets, which improves~\cite{nonlocal-approximation}*{Proposition 3.2}. 
In the nonlocal setting, related iteration techniques were introduced by Ros-Oton and Serra~\cite{ros-oton-serra}, with subsequent extensions by~\cite{fazly}.
The architecture of Proposition~\ref{prop:improve-the-exponent} follows the argument of Proposition~\ref{prop:stable-toda}, but the nonlocal problem carries additional terms, which we prove to have lower-order contributions.
Notably, Lemma~\ref{lemma:interaction-estimate} shows that the contribution from sheet interactions can be absorbed in our estimates.
Combining Proposition~\ref{prop:improve-the-exponent} with the ensuing arguments of~\cite{nonlocal-approximation}, we then convert this integral estimate into pointwise separation of the sheets and establish Theorem~\ref{thm:sheet-happens}.
Accordingly, the optimal seven-dimensional sheeting established in Theorem~\ref{thm:sheet-happens} can be viewed as the nonlocal analogue of the celebrated structure theorem of Schoen and Simon~\cite{sheeting-theorem}.

Let us explain the origin of the $L^{-6}$ bound of Proposition~\ref{prop:improve-the-exponent}.
Schematically, the exponent $p < 4+ 2 \sqrt{2}$ for the $L^{-p}$ estimate of Proposition~\ref{prop:stable-toda} arises from the following observation: for $q := \frac{p-1}{2} > 0$, applying Young's inequality in the iteration argument requires a positive coefficient
\begin{equation}\label{eqn:q-coefficient}
2q + 1 - q^2 > 0, \qquad \iff \qquad q < 1+ \sqrt{2}.
\end{equation}
In the nonlocal setting, the difference between successive graphical sheets, normalized by the natural scale $\sqrt{1-s}$, satisfies properties analogous to~\eqref{eqn:weak-supersolution-semilinear} - \eqref{eqn:stability-semilinear}.
Thus, our argument proceeds by a nonlocal version of this iteration, aimed at obtaining an $L^{-2q-1-s}$ estimate for the function $w := \frac{g_{i+1} - g_i}{\sqrt{1-s}}$.
In particular,~\cite{nonlocal-approximation}*{Proposition 3.2} can be viewed as the $q=0$ case of this result.
In our situation, the corresponding term appearing in~\eqref{eqn:q-delta-epsilon} is a small perturbation of the coefficient~\eqref{eqn:q-coefficient}, provided that certain additional parameters are taken to be suitably small and $s$ is taken sufficiently close to $1$.
In particular, this allows us to set $q = \frac{21}{10}$ and make $2q+1+s>6$.
%, thereby obtaining an $L^{-6}$ estimate.

We then turn to the structure of stable $s$-minimal cones in $\mathbb{R}^n$, for $n \leq 8$.
By the stable Bernstein theorem for minimal surfaces in $\mathbb{R}^6$ and $\mathbb{R}^7$~\cites{mazet , hong-li-wang }, Theorem~\ref{thm:sheet-happens} allows us to apply the regularity properties of~\cite{nonlocal-approximation} and study the spherical links of the cones.
The links of non-flat cones in a contradiction sequence would converge to a minimal lamination $\cL$ of $\bS^{N-1}$, by~\cite{meeks-perez-ros}.
The two-sided cover of a limit leaf of $\cL$ would be stable, contradicting~\cites{catino-mastrolia-roncoroni , wang-bernstein } for $n -2 \leq 6$.
Then, $\cL$ cannot have a limit leaf, so it consists of one compact minimal hypersurface whose cone is stable in $\bR^n$.
By standard theory, for $n \leq 7$, the limit leaf $\cL$ must be an equator~\cite{simons}.

Finally, with these results in place, we can proceed as in the original argument of~\cite{nonlocal-approximation}*{\S 4}, using the sharp Hardy inequality for the limiting problem.
We carry out this step in Section~\ref{sec:completion}.
This final part of the argument was already observed in~\cite{nonlocal-approximation} to apply in all dimensions $n \leq 7$, provided that all other estimates of Theorem~\ref{thm:sheet-happens} were available.
In fact, these arguments extend to $\mathbb{R}^8$, leading to the compactness Theorem~\ref{thm:eight-dimensions} for stable $s$-minimal cones in $\mathbb{R}^8$.
The additional ingredient in eight dimensions is the observation that the ``numerology'' of the limiting system forces either multiplicity-one convergence to the link of a stable classical cone, or multiplicity-two convergence together with a Jacobi field subject to a limiting stability condition.
In the latter case, testing stability with $\phi = |A_{\Sigma}|$ and using Kato's inequality, as in the original argument of Simons~\cite{simons}, proves $|A_{\Sigma}| \equiv 0$ and forces the multiplicity-two limit to be a flat hyperplane.

\smallskip \noindent \textbf{Acknowledgments}.
We are thankful to Hardy Chan for many helpful discussions and comments, as well as for his contribution to the development of nonlocal approximation theory.

SD was supported by the Australian Research Council Future Fellowship FT230100333
``New perspectives on nonlocal equations''. 
JS was supported by the European Research Council under the Grant Agreement No 948029.
RT was supported by a Simons Dissertation Fellowship, an A.G.~Leventis Foundation Scholarship, and an Onassis Foundation Scholarship.
EV was supported by the Australian Laureate Fellowship FL190100081 ``Minimal
surfaces, free boundaries and partial differential equations''.

\smallskip \noindent \textbf{AI statement}.
No stage of this work involved LLMs. 
A separate literature search (Google Gemini), queried for works carrying out an iteration similar to the one of Meadows~\cite{meadows} in the nonlocal setting, led to the reference~\cite{fazly}. The ideas of~\cite{fazly}*{\S 4} helped design the iteration of Proposition~\ref{prop:improve-the-exponent}, whose proof ultimately relies on the ideas of~\cite{ros-oton-serra} and~\cite{nonlocal-approximation}.

\section{Preliminaries from nonlocal approximation theory}\label{sec:prelimiinaries-nonlocal-approximation}

We recall some important properties from~\cite{nonlocal-approximation}, which will be useful in establishing Proposition~\ref{prop:improve-the-exponent}.
We adopt the graphical setting and results developed therein, using only the notation and quantities that are relevant for the proof.

We write points of $\bR^n$ as $x = (x', x_n) \in \bR^{n-1} \times \bR$ and consider $(n-1)$-dimensional balls $B_r^{n-1} := \{ |x'| < r \} \subset \bR^{n-1}$.
Let $E \subset \bR^n$ be a stable $s$-minimal set in a fixed cylinder $\Omega := B_{1/2}^{n-1} \times (-1,1) \subset \bR^n$ whose boundary $\partial E$ in $\Omega$ coincides with the union of graphs $\{ \Gamma_j \}_{i=1}^N$, namely
\[
\partial E \cap \Omega = \bigcup_{j=1}^N \Gamma_j,\qquad{\mbox{with}} \qquad \Gamma_j :=  \{ (x', g_j(x')) : x' \in B_{1/2}^{n-1}, \; -1<g_j(x') < 1 \}
\]
inside $\Omega$, for functions $g_i : B^{n-1}_1 \to \bR$ satisfying
\begin{equation}\label{eqn:slope-bound-g-i}
g_1<g_2<\cdots< g_N\qquad{\mbox{and}} \qquad \sup_j ( \|D g_j \|_{L^{\infty}} + \|D^2 g_j \|_{L^{\infty}}) \leq \delta.
\end{equation}
Here, $N$ is an arbitrarily large positive integer and  $\delta>0$ is sufficiently small.

We express the fractional $s$-mean curvature of a Caccioppoli set $E \subseteq \bR^n$ with the convention
\[
H_s[E](x)= \textup{PV}\int_{\mathbb{R}^n} \frac{\mathbf{1}_{E^c}(y)-\mathbf{1}_E(y)}{|x-y|^{n+s}}\,dy.
\]
Let us record a consequence of the gradient and Hessian bounds in~\eqref{eqn:slope-bound-g-i}, cf.~\cite{nonlocal-approximation}*{Lemma 3.3}.
\begin{lemma}\label{lemma:gradient-hessian-bound}
Inside the ball $B_{1/8}^{n-1}(0)$, we have the bound
\begin{equation}\label{eqn:gj-x-prime}
    |D (g_j - g_i)| \leq 8 \delta^{\frac{1}{2}} (g_j - g_i)^{\frac{1}{2}}.
\end{equation}
\end{lemma}
\begin{proof}
Let $x' \in B^{n-1}_{1/8}(0)$. We consider the function 
\[
h(z') := g_j(x'+z') - g_i(x'+z').
\]
The assumptions on the gradient and the Hessian of the $g_j$ yield $|Dh| \leq \delta$ and $|D^2 h| \leq \delta$, therefore
\[
0 \leq h(z') = h(0) + \la Dh(0), z' \rg + \delta \, |z'|^2 \qquad \text{for } \; |z'| \leq \tfrac{1}{4}.
\]
For every $x' \in B^{n-1}_{1/8}(0)$ and every $r \in (0, \frac{1}{4})$, the ball $B^{n-1}_r(x')$ is contained in the domain $B^{n-1}_{1/2}(0)$ of the functions $g_i$, hence $B^{n-1}_r(0)$ is contained in the domain of $h$.
We therefore obtain, for $r \in (0, \frac{1}{4})$,
\[
|Dh(0)| = r^{-1} \sup_{z' \in B^{n-1}_r(0)} ( - \la Dh(0), z' \rg) \leq r^{-1} h(0) + \delta r.
\]
If $h(0) < \frac{1}{16} \delta$, then we can choose $r = \sqrt{\delta^{-1} h(0)} \in (0, \frac{1}{4})$ in this inequality to conclude that $|Dh(0)| \leq 2 \sqrt{h(0) \delta}$.
We therefore obtain
\[ |D(g_j - g_i)(x')| \leq 2 \delta^{\frac{1}{2}} \sqrt{(g_j - g_i)(x')} \, .
\]
In addition, if $h(0) \geq \frac{1}{16} \delta$, then the bound $|Dh| \leq \delta$ implies that
\[
|Dh(0)| \leq 2 \delta \leq 8 \textstyle{\sqrt{\delta h(0)}}.
\] 
Therefore, using the inequality $\sqrt{\delta \, h(0)} \geq \frac{1}{4} \delta$, we conclude that
\[ 
|D(g_j - g_i)(x')| \leq 8 \delta^{\frac{1}{2}} \textstyle{\sqrt{(g_j - g_i)(x')}}.
\]
The claimed bound~\eqref{lemma:gradient-hessian-bound} follows from these inequalities.
\end{proof}

Given a domain $D$ and a measurable symmetric kernel $K$ on $D \times D$, we define $\cL_K$ weakly by
\begin{equation}\label{eqn:Lu-kernel}
    \int_D v \cL_K u = \frac{1}{2} \iint_{D \times D} (u (x') - u(y')) \, (v(x') - v(y')) \, K(x',y') \, dx' \, dy'.
\end{equation}
For the kernels considered below and $u \in C^2(D)$, $\cL_K$ admits the pointwise representation
\begin{equation}\label{eqn:LK-operator-over-D}
\cL_K u(x') = \textup{PV} \int_D (u(x') - u(y')) \, K(x', y') \, dy'.
\end{equation}

In what follows, we pick $i \in \{ 1, \dots, N-1\}$ such that $|g_i|, |g_{i+1}| \leq \frac{1}{2}$ on $B^{n-1}_{1/2}$ and form the graphs over $B^{n-1}_{1/2}$, given by
\[
\Gamma_j = \{ (x', g_j(x')) : x' \in B^{n-1}_{1/2}, \; -1<g_j(x') < 1\}.
\]
We also let $c_{n,s} := \int_{\bR^{n-1}} ( 1 + |x'|^2)^{- \frac{n+s}{2}} \, dx'$.
As in~\cite{nonlocal-approximation}*{\S 3.2}, the next step towards the proof of Proposition~\ref{prop:improve-the-exponent} will be to localize the stability inequality to the interaction between two successive sheets $\Gamma_i, \Gamma_{i+1}$ of $E$.
First, we will need a one-dimensional kernel bound in the style of Riesz's rising sun lemma; see, for example,~\cite{anisotropic-michael-simon}*{Lemma 3.11}.
\begin{lemma}\label{lemma:rising-sun-lemma}
Let $k \in C^1(\bR)$ be a non-negative, even, non-increasing function on $(0,\infty)$ with $k(t) \to 0$ as $|t| \to \infty$.
If $f$ takes values in $[0,1]$ and vanishes on $(a,b)$, then
\[
- \int_{\bR} f(t) k'(t) \, dt \leq k(b) + k(0) - k(a).
\]
\end{lemma}
\begin{proof}
Since $-k'(t)\leq0$ for $t<0$, we may discard the contribution from $(-\infty,0)$. 
On $(0,\infty)$, $d\mu(t) :=-k'(t),dt$ is a positive measure. 
Since $f$ takes values in $[0,1]$ and vanishes on $(a,b)$, the result follows immediately by considering the three cases $b<0$, $0 \in (a,b)$, and $a>0$.
\end{proof}

Next, we recall an estimate for the contribution from exterior terms as in~\cite{nonlocal-approximation}*{Appendix A}.

\begin{lemma}\label{lemma:slab-estimate}
Let $S$ be the slab enclosed between two successive sheets, namely
\[
S := \{ (x', t) : x' \in B^{n-1}_{1/2}, \; g_i(x') < t < g_{i+1}(x') \}
\]
and suppose that $S \subset E$.
Then, 
\begin{enumerate}[(i)]
    \item For all $x = (x', x_n) \in (B^{n-1}_{1/4} \times \bR) \cap \Gamma_i$, we have that
\begin{equation}\label{pr:tw:01}
\Bigl \la \nu(x), \int_{E \setminus S} \nabla_y |x-y|^{-n-s} \, dy \Bigr \rg \leq (1+C \delta) \, c_{n,s}(g_{i+2} - g_i)(x')^{-1-s} + C \Bigl( 1 + \frac{\delta^2}{1-s} \Bigr)
\end{equation}
whenever $i \leq N-2$.
\item Likewise, if $x = (x', x_n) \in (B^{n-1}_{1/4} \times \bR) \cap \Gamma_{i+1}$, then
\begin{equation}\label{pr:tw:02}
\Bigl\la \nu(x), \int_{E \setminus S} \nabla_y |x-y|^{-n-s} \, dy \Bigr \rg \leq (1+C \delta) \, c_{n,s} \, (g_{i+1} - g_{i-1})(x')^{-1-s} + C \Bigl( 1 + \frac{\delta^2}{1-s} \Bigr)
\end{equation}
with the first term omitted when $i=1$.
\end{enumerate}
\end{lemma}
\begin{proof}
This result is an immediate consequence of Step 1 of the proof of~\cite{nonlocal-approximation}*{Proposition 3.5}.
It suffices to record the estimate in~\eqref{pr:tw:01}, since the one in~\eqref{pr:tw:02} then follows after reversing the vertical direction.
We fix $x \in \Gamma_i$ with $x' \in B^{n-1}_{1/4}$ and apply a rigid motion to assume that $x =0$, that the tangent hyperplane to $\Gamma_i$ is $T_x \Gamma_i = \{ x_n = 0 \}$, and that the normal vector is $\nu(x) = - e_n$.

In a given cylinder $B^{n-1}_{r_0} \times (-r_0,r_0)$, the sheets remain ordered graphs $\tilde{g}_0<\tilde{g}_1 < \tilde{g}_2 < \cdots$, with $\tilde{g}_k$ representing $\Gamma_{i+k}$ for $k \geq 0$, and 
\begin{equation}\label{eqn:localize-h0}
\tilde{g}_0(0) = 0, \qquad D\tilde{g}_0(0) = 0, \qquad |\tilde{g}_0(y')| \leq C \delta |y'|^2.
\end{equation}
The contribution to the integral from outside this cylinder is bounded by a uniform constant $C$; for the region inside the cylinder, we use $\nu(x) = -e_n$ to compute
\begin{equation}\label{eqn:nu(x)-W-r}
\la \nu(x), \nabla_y |y|^{-n-s} \rg = - W'_{|y'|}(t), \qquad \text{where } \quad W_{|y'|}(t) = (|y'|^2 + t^2)^{- \frac{n+s}{2}}.
\end{equation}
Here, $\textup{sgn} \, W'_r(t) = - \textup{sgn}(t)$.
Because $S \subset E$, the region enclosed between $\tilde{g}_0$ and $\tilde{g}_1$ is absent from $E \setminus S$.
Along every vertical fiber, the remaining set can therefore contribute positively only 
\begin{enumerate}[(i)]
    \item between $t=0$ and $t=\tilde{g}_0(y')$, when $\tilde{g}_0(y')>0$; and
    \item above the first sheet beyond the selected pair, namely $\tilde{g}_2(y')$, when that sheet exists.
\end{enumerate}
Let $r = |y'|$.
Combining these computations, we can use Lemma~\ref{lemma:rising-sun-lemma} with $k(t) = W_r(t)$ to obtain
\begin{equation}\label{eqn:rising-sun-lemma}
\int_{ \{ t : (y', t) \in E \setminus S\} } - W'_r(t) \, dt \leq W_r (\tilde{g}_2(y')) + ( W_r(0) - W_r(\tilde{g}_0(y'))
\end{equation}
with the first summand $W_r(\tilde{g}_2(y'))$ omitted if $i=N-1$.
Notice that~\eqref{eqn:rising-sun-lemma} continues to hold if $\tilde{g}_2(y')<0$, since $\tilde{g}_0(y') < \tilde{g}_2(y')<0$ makes
\[
W_r(0) - W_r(\tilde{g}_2(y')) \leq W_r(0) - W_r(\tilde{g}_0(y')).
\]
To estimate the terms in~\eqref{eqn:rising-sun-lemma}, we use the bound~\eqref{eqn:localize-h0} and the evenness of $W_r$ to conclude that
\begin{equation}\label{eqn:W-integral-estimates}
    \begin{split}
W_r(0) - W_r(\tilde{g}_0(y')) &\leq C \, \tilde{g}_0(y')^2 \, r^{-n-s-2} \leq C \, \delta^2 r^{2-n-s} \\ {\mbox{and }}\quad
\int_{B_{r_0}^{n-1}} (W_r(0) - W_r(\tilde{g}_0(y')) \, dy' &\leq C \delta^2 \int_0^{r_0} r^{-s} \, dr \leq \frac{C \delta^2}{1-s}.
    \end{split}
\end{equation}
Moreover, the graph $\tilde{g}_2$ has slope $O(\delta)$, hence we can use the estimate of~\cite{nonlocal-approximation}*{Lemma 3.3} to control the first term in the right-hand side of~\eqref{eqn:rising-sun-lemma} as
\begin{equation}\label{eqn:W-r-h2(y')}
\int_{B_{r_0}^{n-1}} W_r (\tilde{g}_2(y')) \, dy' \leq (1+C \delta) \, c_{n,s} \, (g_{i+2} - g_i)(x')^{-1-s} + C.
\end{equation}
Here, we included an additive dimensional constant from truncating the graph to the fixed cylinder.
Combining the bounds~\eqref{eqn:nu(x)-W-r}-\eqref{eqn:W-r-h2(y')} and adding the bounded contribution outside the cylinder, we find
\[
\Bigl\la \nu(x), \int_{E \setminus S} \nabla_y |x-y|^{-n-s} \, dy \Bigr\rg \leq (1+C\delta) \, c_{n,s} \, (g_{i+2} - g_i)(x')^{-1-s} + C \Bigl( 1 + \frac{\delta^2}{1-s} \Bigr).
\]
If $i=N-1$, there is no sheet $\Gamma_{i+2}$, and the corresponding term is absent.
For $x \in \Gamma_{i+1}$, we reverse the vertical direction, so the first sheet beyond the selected pair is $\Gamma_{i-1}$, and the same argument gives the second estimate.
This completes the proof of Lemma~\ref{lemma:slab-estimate}.
\end{proof}

The following result is a consequence of Sections~3.1, 3.2, and 3.3 of~\cite{nonlocal-approximation}.

\begin{proposition}\label{prop:first-prop}
There are dimensional constants $\delta_* > 0$ and $s_*\in( 0,1)$ with the following property.
For $\delta \in (0,\delta_*)$ and $s \in (s_*, 1)$ in our setting, there  is a symmetric kernel $K_i$ on $B^{n-1}_{1/2} \times B^{n-1}_{1/2}$ with
\begin{equation}\label{eqn:K-kernel-property}
(1- C\delta^2) \, |x'-y'|^{-n-s} \leq K_i (x', y') \leq |x' - y'|^{-n-s}
\end{equation}
such that the associated operator $\cL_{K_i}$ satisfies
\begin{equation}\label{eqn:L-kernel-property}
- \cL_{K_i} (g_{i+1} - g_i) \leq ( 2 c_{n,s} s^{-1} + C \delta) \, |g_{i+1} - g_i|^{-s} + C \qquad \text{on } \; B_{1/4}^{n-1}.
\end{equation}
Moreover, for every $\phi \in C_c^1(B^{n-1}_{1/4})$, it holds that
\begin{equation}\label{eqn:stability-iteration}
    \begin{split}
&(2-C \delta) \, c_{n,s} \int_{B^{n-1}_{1/2}} \frac{\phi^2}{|g_{i+1} - g_i|^{1+s}} \\
&\leq (1+C \delta^2) \int_{B^{n-1}_{1/2}} \phi \cL_{K_i} \phi + \frac{1}{2} \iint_{\Gamma_i \times \Gamma_{i+1}} \frac{( \phi(x') - \phi(y'))^2}{|x-y|^{n+s}} \, d \cH^{n-1}_x \, d \cH^{n-1}_y \\
& \quad + \Bigl( \frac{1}{2} + C \delta \Bigr) \, c_{n,s} \int_{B^{n-1}_{1/2}} \left[ \frac{1}{(g_{i+1} - g_{i-1})^{1+s}} + \frac{1}{(g_{i+2} - g_i)^{1+s}} \right] \phi^2 + C \Bigl( 1 + \frac{\delta^2}{1-s} \Bigr) \int_{B^{n-1}_{1/2}} \phi^2.
    \end{split}
\end{equation}
The constant $C$ depends only on $n$.
\end{proposition}
\begin{proof}
We proceed in two steps: first, we use~\cite{nonlocal-approximation}*{Proposition 3.5} to produce the kernel $K_i$ and obtain the estimates~\eqref{eqn:K-kernel-property} and~\eqref{eqn:L-kernel-property}.
Then, we apply~\cite{nonlocal-approximation}*{Proposition 3.1 and Lemma 3.3} to obtain the stability inequality~\eqref{eqn:stability-iteration}.

\smallskip \noindent \textbf{Step 1:}
For two adjacent sheets $g_i, g_{i+1}$, we consider the associated mean value kernel
\[
K_i(x',y') := \int_0^1 \bigl( |x'-y'|^2 + [ (1-t) (g_i(y') - g_i(x')) + t \, ( g_{i+1}(y') - g_{i+1}(x')) ]^2 \bigr)^{- \frac{n+s}{2}} \, dt
\]
which arises from the graphical linearization of the $s$-mean curvature operator, obtained in~\cite{nonlocal-approximation}*{(3.34)--(3.36)}.
The symmetry of $K_i(x',y')$ is automatic, and the slope bound~\eqref{eqn:slope-bound-g-i} on the $g_i$ yields
\[
|(1-t) \, (g_i(y') - g_i(x')) + t \, (g_{i+1}(y') - g_{i+1}(x'))| \leq \delta \, |x'-y'|.
\]
We therefore obtain $(1+\delta^2)^{- \frac{n+s}{2}} |x'-y'|^{-n-s} \leq K_i(x', y') \leq |x'-y'|^{-n-s}$, and combining this inequality with the bound $(1+\delta^2)^{- \frac{n+s}{2}} \geq 1 - C_n \delta^2$ yields the first property in~\eqref{eqn:K-kernel-property}.

Next, replacing $E$ by its complement $E^c$, if necessary, we can assume that the slab
\[
S := \{ (x', t) : x' \in B^{n-1}_{1/2}, \; g_i(x') < t < g_{i+1}(x') \}
\]
enclosed between the two graphs is contained in $E$.
We pick a point $x'$ and set
\[
x_- := (x', g_i(x')), \qquad x_+ := (x', g_{i+1}(x')).
\]
Then, since $S \subset E$ and $E$ is $s$-stationary, the monotonicity of $H_s$ under containment implies that 
\[
H_s[S](x) - H_s[E](x) = 2 \int_{E \setminus S} |x-y|^{-(n+s)} \, dy \geq 0
\]
and hence $H_s[S](x_{\pm}) \geq 0$.
We also recall that
\[
g_{i+1}(y') - g_i(x') = (g_{i+1} - g_i)(x') + ( g_{i+1}(y') - g_{i+1}(x')) \geq (g_{i+1} - g_i) (x') - \delta \, |x'-y'|,
\]
due to the gradient bound $|Dg_{i+1}| \leq \delta$.
Combining this property with the fact that $H_s[S](x_{\pm}) \geq 0$, and with~\cite{nonlocal-approximation}*{(3.34)-(3.35)}, we see that, for all $x' \in B^{n-1}_{1/4}$,
\begin{equation}\label{eqn:integral-kernel-bound}
- \cL_{K_i} (g_{i+1} - g_i) \leq 2 \int_{\bR^{n-1}} \int^{\infty}_{g_{i+1}(x') - g_i(x') - \delta |y'|} \frac{dt \, dy'}{( |y'|^2 + t^2)^{\frac{n+s}{2}}} + C.
\end{equation}

The symmetric argument shows that
\[
g_{i+1}(x') - g_i(y') \geq (g_{i+1} -g_i)(x') - \delta \, |x'-y'|.
\]

Now, to estimate the right-hand side of the bound~\eqref{eqn:integral-kernel-bound}, we claim that
\begin{align}
\int_{\bR^{n-1}} \int_a^{\infty} ( |y'|^2 + t^2)^{- \frac{n+s}{2}} \, dt \, dy' &= c_{n,s} \, s^{-1} \, a^{-s}\label{eqn:a-integral-equality} \end{align} and
\begin{align}
\int_{\bR^{n-1}} \int^1_{1 - \delta |z|} (|z|^2+t^2)^{- \frac{n+s}{2}} \, dt \, dz &\leq C_n \delta, \label{eqn:a-tail-inequality}
\end{align}
for a constant $C_n$ depending only on $n$, uniformly in $s  \in [\frac{1}{2},1)$ and $\delta \in (0,\frac{1}{4}]$.

The equality~\eqref{eqn:a-integral-equality} follows by integration. 
For the bound in~\eqref{eqn:a-tail-inequality}, we split the integral into $|z| \leq 2$ and $|z|>2$.
If $|z| \leq 2$, then $t \geq 1-2 \delta \geq \frac{1}{2}$, and the interval of integration has length $\delta |z|$, hence
\[
\int_{|z| \leq 2} \int_{1 - \delta |z|}^1 \frac{dt \, dz}{(|z|^2+t^2)^{\frac{n+s}{2}}} \leq 2^{n+1} \delta \int_{|z| \leq 2} \, |z| \, dz = \frac{2^{2n+1}}{n} |\bS^{n-2}| \delta
\]
by $n+s \leq n+1$.

If instead $|z|>2$, we use $(|z|^2+t^2){- \frac{n+s}{2}} \leq |z|^{-n-s}$ on an interval of length $\delta |z|$ to bound
\begin{align*}
    \int_{|z|>2} \int_{1 - \delta |z|}^1 \frac{dt \, dz}{(|z|^2+t^2)^{\frac{n+s}{2}}} &\leq |\bS^{n-2}| \, \delta \int_2^{\infty} r^{-1-s} \, dr = |\bS^{n-2}| 2^{-s} s^{-1} \delta \leq \sqrt{2} \, |\bS^{n-2}| \delta
\end{align*}
due to $s \in [\frac{1}{2},1)$.
Thus, the second bound~\eqref{eqn:a-tail-inequality} holds with $C_n = \frac{2^{2n+1}}{n} + \sqrt{2}$.

We now combine the estimates~\eqref{eqn:a-integral-equality} and~\eqref{eqn:a-tail-inequality} into the kernel bound~\eqref{eqn:integral-kernel-bound} to obtain
\[
- \cL_{K_i} (g_{i+1} - g_i) \leq ( 2 c_{n,s} s^{-1} + C \delta) ( g_{i+1} - g_i)^{-s} + C.
\]
For $\delta_*>0$ sufficiently small and $s_* < 1$ sufficiently close to $1$, we have $2 c_{n,s} s^{-1} + C \delta \leq \frac{41}{20} c_{n,s}$ whenever $s \in (s_*,1)$ and $\delta \in (0,\delta_*)$.
Consequently,
\begin{equation}\label{eqn:LK-wi-estimate}
    - \cL_{K_i} (g_{i+1} - g_i) \leq \tfrac{41}{20} \, c_{n,s} \, (g_{i+1} - g_i)^{-s} + C.
\end{equation}

\noindent \textbf{Step 2:}
Next, we obtain the bound~\eqref{eqn:stability-iteration} as a byproduct of~\cite{nonlocal-approximation}*{Proposition 3.1}.
Specifically, recalling the notation $\Gamma := \partial E \cap \Omega$ and $\Gamma_i := \{ x_n = g_i(x') \} \cap \Omega_0$ therein, we apply~\cite{nonlocal-approximation}*{Proposition 3.1} with normal speed $\phi(x')$ supported on $\Gamma_i \cup \Gamma_{i+1}$ and equal to zero on the remaining sheets, namely
\[
\eta(x) = \phi(x') \quad \text{for } \; x \in \Gamma_i \cup \Gamma_{i+1}, \qquad \eta(x) = 0 \quad \text{for } \; x \not\in (\Gamma_i \cup \Gamma_{i+1}).
\]
This configuration can clearly be extended to a compactly supported $C^1$ function in $\Omega$.
Separating the interactions for which both variables lie on the two selected sheets, we obtain
\begin{equation}\label{eqn:0-(1-s)-Ext-Omega-2}
0 \leq (1-s) \sum_{\ell,\ell' \in \{ i,i+1\} } \iint_{\Gamma_{\ell} \times \Gamma_{\ell'}} \frac{(\phi(x') - \phi(y'))^2 - |\nu(x) - \nu(y)|^2 \phi(x')^2}{|x-y|^{n+s}} \, d\cH^{n-1}_x \, d \cH^{n-1}_y + \cR_i.
\end{equation}
Here, $\cR_i$ is a remainder term recording mixed interactions as well as exterior terms: using the definition of the term $\textup{Ext}^{\Omega}_2$ from~\cite{nonlocal-approximation}*{Proposition 3.1}, this is given by
\begin{align*}
    \cR_i &= 2(1-s) \Biggl[ \int_{\Gamma_i} \phi(x')^2 \la \nu(x), R(x) \rg \, d \cH^{n-1}_x + \int_{\Gamma_{i+1}} \phi(x')^2 \la \nu(x), R(x) \rg \, d \cH^{n-1}_x \Biggr], 
    \intertext{where we denoted}
    R(x) &:= \int_{\Gamma \setminus (\Gamma_i \cup \Gamma_{i+1})} \frac{\nu(y)}{|x-y|^{n+s}} \, d\cH^{n-1}_y + \int_{E \setminus \Omega} \nabla_y \Bigl( \frac{1}{|x-y|^{n+s}} \Bigr) \, dy + \int_{E \cap \partial \Omega} \frac{n_{\Omega}(y)}{|x-y|^{n+s}} \, d \cH^{n-1}_y.
\end{align*}
Indeed, this term comes from expressing $2 - |\nu(x) - \nu(y)|^2 = 2 \la \nu(x), \nu(y) \rg$ for both orderings of every pair between points of $(\Gamma_i \cup \Gamma_{i+1})$ and of $\Gamma \setminus ( \Gamma_i \cup \Gamma_{i+1})$.
The last two terms in $R(x)$ are precisely those defining $\textup{Ext}_2^{\Omega}$, since 
\[
\nabla_y |x-y|^{-n-s} = (n+s) (x-y) |x-y|^{-(n+s+2)}.
\]
For a given point $x \in \Gamma_i \cup \Gamma_{i+1}$ on a selected sheet, discarding the slab $S$ inside $\Omega$ removes the two selected sheets $\Gamma_i \cup \Gamma_{i+1}$, making the boundary of $E \setminus S$ inside $\Omega$ consist of the unselected sheets with their original outward normals.
For a given point $x$ on a selected sheet, the set $E \setminus S$ is empty in a neighborhood of $x$, and consequently we can apply Green's theorem~\cite{maggi}*{Proposition~19.22 and Lemma~22.11} with no singularity at $x$, obtaining that
\[
\int_{\Gamma \setminus (\Gamma_i \cup \Gamma_{i+1})} \frac{\nu(y)}{|x-y|^{n+s}} \, d \cH^{n-1}_y = \int_{(E \setminus S) \cap \Omega} \nabla_y |x-y|^{-n-s} \, dy - \int_{\partial \Omega} \frac{\mathbf{1}_{E \setminus S}(y) \, n_{\Omega}(y)}{|x-y|^{n+s}} \, d \cH^{n-1}_y.
\]

We now substitute this expression into the definition of $R(x)$, using the fact that $(E \setminus S) \setminus \Omega = E \setminus \Omega$ and $\mathbf{1}_E - \mathbf{1}_{E \setminus S} = \mathbf{1}_S$ due to $S \subset \Omega$.
In this way, we arrive at
\[
\cR_i = 2(1-s) \sum_{\ell=i}^{i+1} \int_{\Gamma_{\ell}} \phi(x')^2 \, \left\la \nu(x), \Biggl[ \int_{E \setminus S} \nabla_y |x-y|^{-n-s} \, dy + \int_{\partial \Omega} \frac{\mathbf{1}_S(y) \, n_{\Omega}(y)}{|x-y|^{n+s}} \, d \cH^{n-1}_y \Biggr] \right\rg \, d \cH^{n-1}_x
\]
In the next step, we notice that the trace of $S$ is supported on the lateral boundary of $\Omega$, whose distance from a selected point with $\phi(x') \neq 0$ is at least $\frac{1}{4}$.
Thus, for a uniform constant~$C_n$, by arguing as in~\cite{nonlocal-approximation}*{Appendix A} we see that
\begin{align*}
    \left| \int_{\partial \Omega} \frac{\mathbf{1}_S(y) \, n_{\Omega}(y)}{|x-y|^{n+s}} \, d \cH^{n-1}_y \right| + \int_{(E \setminus S) \setminus \Omega} \bigl| \nabla_y |x-y|^{-n-s} \bigr| \, dy &\leq C_n 
\end{align*}

In addition, the assumption in Step 1 implies that the slab region $S$ enclosed between the two selected adjacent graphs is contained in $E$, and accordingly the remaining terms over $\Omega$ are estimated by Lemma~\ref{lemma:slab-estimate}.
Furthermore, we notice that the two selected sheets $\Gamma_i, \Gamma_{i+1}$ can be parametrized over $B^{n-1}_{1/4}$ with surface Jacobians of order $1+O(\delta^2)$, therefore $ d\cH^{n-1}_x =(1 + O(\delta^2)) \, dx'$ on the selected sheets.
Thus, combining these estimates, we conclude that
\begin{equation}\label{eqn:for-final-I-first}
    \begin{split}
    \cR_i &\leq (2+C \delta) \, c_{n,s}\, (1-s) \int_{B^{n-1}_{1/4}} [ (g_{i+2} - g_i)^{-1-s} + (g_{i+1} - g_{i-1})^{-1-s} ] \, \phi^2 \, dx' \\
    & \quad + C (\delta^2+1-s) \int_{B^{n-1}_{1/4}} \phi^2 \, dx'.
    \end{split}
\end{equation}
The term with $g_{i-1}$ is present when $i \geq 2$, and the term with $g_{i+2}$ is present when $i \leq N-2$.

It remains to estimate the terms in~\eqref{eqn:0-(1-s)-Ext-Omega-2} involving only the selected sheets $\Gamma_i$ and $\Gamma_{i+1}$.
Discarding the non-positive normal term $- |\nu(x) - \nu(y)|^2$ and applying the gradient bounds $|Dg_j| \leq \delta$ for the graphs $\Gamma_j$, as well as the kernel bound~\eqref{eqn:K-kernel-property}, we can bound the self-interaction part of~\eqref{eqn:0-(1-s)-Ext-Omega-2} as
\begin{align*}
    & \sum_{j=i}^{i+1} \iint_{\Gamma_j \times \Gamma_j} \frac{(\phi(x') - \phi(y'))^2 - |\nu(x) - \nu(y)|^2 \phi(x')^2}{|x-y|^{n+s}} \, d\cH^{n-1}_x \, d \cH^{n-1}_y \\
    &\leq (2+C\delta^2)\iint_{B^{n-1}_{1/2} \times B^{n-1}_{1/2}} ( \phi(x') - \phi(y'))^2 \, K_i(x',y') \, dx' \, dy'.
\end{align*}
Here, we used the graphical expression of $\Gamma_j$ over $\Gamma_j$ in addition to the bound~\eqref{eqn:K-kernel-property} for the kernel $K_i$.
In view of the definition~\eqref{eqn:Lu-kernel}, the right-hand side is expressible in terms of the kernel $\cL_{K_i}$ as
\begin{equation}\label{eqn:for-final-I}
    \begin{split}
    \sum_{j=i}^{i+1} \iint_{\Gamma_j \times \Gamma_j} \frac{(\phi(x') - \phi(y'))^2 - |\nu(x) - \nu(y)|^2 \phi(x')^2}{|x-y|^{n+s}} \, d\cH^{n-1}_x \, d \cH^{n-1}_y \leq 4 (1+C \delta^2) \int_{B^{n-1}_{1/2}} \phi \cL_{K_i} \phi.
    \end{split}
\end{equation}
Thus, in the next step, we estimate the interaction between the two selected sheets $\Gamma_i$ and $\Gamma_{i+1}$.
Because the above stability inequality involves an integral over ordered pairs, the combined contribution from $\Gamma_i \times \Gamma_{i+1}$ and $\Gamma_{i+1} \times \Gamma_i$ to the first term appearing in~\eqref{eqn:0-(1-s)-Ext-Omega-2} is
\begin{equation}\label{eqn:interacting-successive}
    \begin{split}
    & \sum_{ \{ \ell, \ell' \} = \{ i,i+1\} } \frac{(\phi(x') - \phi(y'))^2}{|x-y|^{n+s}} \, d \cH^{n-1}_x \, d \cH^{n-1}_y  = 2 \iint_{\Gamma_i \times \Gamma_{i+1}} \frac{(\phi(x') - \phi(y'))^2}{|x-y|^{n+s}} \, d \cH^{n-1}_x \, d \cH^{n-1}_y.
    \end{split}
\end{equation}
For the normal terms in~\eqref{eqn:0-(1-s)-Ext-Omega-2}, the two orderings of the graphs produce
\begin{equation}\label{eqn:normal-terms-phi(x')-phi(y')}
    \begin{split}
    & \sum_\{ \ell, \ell' \} = \{ i,i+1\} \frac{- |\nu(x) - \nu(y)|^2 \phi(x')^2}{|x-y|^{n+s}} \, d \cH^{n-1}_x \, d \cH^{n-1}_y \\
    &= - \iint_{\Gamma_i \times \Gamma_{i+1}} \frac{|\nu(x) - \nu(y)|^2 \, (\phi(x')^2 + \phi(y')^2)}{|x-y|^{n+s}} \, d \cH^{n-1}_x \, d \cH^{n-1}_y.
    \end{split}
\end{equation}
By the assumption in Step 1, the slab region $S$ enclosed between the two selected adjacent graphs is contained in $E$, so their outward normals point in opposite directions.
Combined with the gradient bound $|Dg_j| \leq \delta$ on the graphs, this yields that, for all~$(x,y) \in \Gamma_i \times \Gamma_{i+1}$,
\begin{equation}\label{eqn:normal-terms-estimate}
|\nu(x) - \nu(y)|^2 = 4 + O(\delta^2) .
\end{equation}
Moreover,~\cite{nonlocal-approximation}*{Lemma 3.3} shows that, for a point $x$ on either selected sheet, e.g., $x \in \Gamma_i$,
\begin{equation}\label{eqn:x-in-amma-i}
\int_{\Gamma_{i+1}} |x-y|^{-(n+s)} \, d \cH^{n-1}_y \geq (1 - C \delta) \, c_{n,s} \, (g_{i+1} - g_i)(x')^{-1-s} - C.
\end{equation}
The bound~\eqref{eqn:x-in-amma-i} is entirely analogous for $x \in \Gamma_{i+1}$, with the integral taken over $y \in \Gamma_i$.

Therefore, by~\eqref{eqn:x-in-amma-i}, each of the two terms in $\phi(x')^2 + \phi(y')^2$ appearing in~\eqref{eqn:normal-terms-phi(x')-phi(y')} produces
\allowdisplaybreaks{
\begin{align*}
    \iint_{\Gamma_i \times \Gamma_{i+1}} \phi(x')^2 \, |x-y|^{-(n+s)} \, d \cH^{n-1}_x \, d \cH^{n-1}_y &\geq (1 - C \delta) \, c_{n,s} \int_{B^{n-1}_{1/2}} (g_{i+1} - g_i)^{-1-s} \phi^2, \\
    \iint_{\Gamma_i \times \Gamma_{i+1}} \phi(y')^2 |x-y|^{-(n+s)} \, d \cH^{n-1}_x \, d \cH^{n-1}_y &\geq (1-C \delta) \, c_{n,s} \int_{B^{n-1}_{1/2}} (g_{i+1} - g_i)^{-1-s} \phi^2.
\end{align*}}
In view of~\eqref{eqn:normal-terms-estimate}, this contribution is then multiplied by a factor of $4 + O(\delta^2)$.
Consequently, the right-hand of the equality~\eqref{eqn:normal-terms-phi(x')-phi(y')} is bounded from above by
\allowdisplaybreaks{
\begin{equation}\label{eqn:for-final-II}
    \begin{split}
& - (1-s) \iint_{\Gamma_i \times \Gamma_{i+1}} \frac{|\nu(x) - \nu(y)|^2 \, (\phi(x')^2 + \phi(y')^2)}{|x-y|^{n+s}} \, d \cH^{n-1}_x \, d \cH^{n-1}_y \\
&\leq - (8 - C \delta) \, c_{n,s} \, (1-s) \int_{B^{n-1}_{1/2}} (g_{i+1} - g_i)^{-1-s} \phi^2 \, dx' + C(1-s) \int_{B^{n-1}_{1/2}} \phi^2 \, dx'.
    \end{split}
\end{equation}}
Finally, we combine~\eqref{eqn:for-final-I-first}, \eqref{eqn:interacting-successive}, and~\eqref{eqn:for-final-II} into~\eqref{eqn:0-(1-s)-Ext-Omega-2}, rearrange the negative term from~\eqref{eqn:for-final-II} to the left, and divide by $4(1-s)$.
We then add the bound~\eqref{eqn:for-final-I}.
This produces the desired inequality~\eqref{eqn:stability-iteration}.
\end{proof}

We conclude this section with an $L^{-p}$ estimate and a linear lower bound for stable weak supersolutions of the equation $\Delta_{\Sigma} w \leq 4 w^{-1}$ on minimal hypersurfaces.
We obtain this result by arguing in the style of Schoen-Simon-Yau~\cite{schoen-simon-yau}.
A similar iteration was carried out by Meadows~\cite{meadows}*{\S 7 - 9} in the context of stable solutions of the problem $\Delta u = \frac{1}{u}$.
Our proof of Proposition~\ref{prop:improve-the-exponent} will be modeled on applying this technique in the nonlocal setting.

\begin{proposition}\label{prop:stable-toda}
Consider a complete hypersurface $\Sigma \subset \bR^n$.
Given a point $x_0 \in \Sigma$ and $R>0$, let $w \in H^1_{\textup{loc}}(B^{\Sigma}_{2R}(x_0))$ be a positive, stable weak supersolution of the equation $\Delta_{\Sigma} w \leq 4 w^{-1}$, namely
\begin{align}
    - \int_{\Sigma} \la \nabla^{\Sigma} w, \nabla^{\Sigma} \phi \rg &\leq 4 \int_{\Sigma} w^{-1} \phi, \qquad & \text{for all } \; \phi & \in C_c^1(B^{\Sigma}_{2R}(x_0)) \text{ with $\phi \geq 0$}, \label{eqn:weak-supersolution-semilinear} \\
    4 \int_{\Sigma} w^{-2} \phi^2 &\leq \int_{\Sigma} |\nabla^{\Sigma} \phi|^2, \qquad & \text{for all } \; \phi & \in C_c^1(B^{\Sigma}_{2R}(x_0)). \label{eqn:stability-semilinear}
\end{align}
Then, for every $p \in (2, 4+2 \sqrt{2})$, there is a universal constant $C_p$ such that
\begin{equation}\label{eqn:L-p-estimate}
\int_{B^{\Sigma}_R(x_0)} w^{-p} \leq C_p \, R^{-p} \, \textup{Vol}_{\Sigma} (B^{\Sigma}_{2R}(x_0)).
\end{equation}
Moreover, when $n \leq 7$ and $\Sigma^{n-1} \subset \mathbb{R}^n$ is a complete minimal hypersurface with Euclidean volume growth, so that $\textup{Vol}_{\Sigma} (B^{\Sigma}_{2R}(x_0)) \leq \Lambda R^{n-1}$, we have $\essinf_{B^{\Sigma}_R(x_0)} w \geq c(n,\Lambda) R$ for $R \geq R_0$.
\end{proposition}
\begin{proof}
For a cutoff function $\eta \in C_c^1(B^{\Sigma}_{2R}(x_0))$, testing the stability inequality with $\phi = w^{-q} \eta$ gives
\begin{align*}
    4 \int_{\Sigma} w^{-2(q+1)} \eta^2 &\leq q^2 \int_{\Sigma} \eta^2 w^{-2(q+1)} |\nabla^{\Sigma} w|^2 + \int_{\Sigma} w^{-2q} |\nabla^{\Sigma} \eta|^2 - 2q \int_{\Sigma} w^{-(2q+1)} \eta \la \nabla^{\Sigma} w, \nabla^{\Sigma} \eta \rg.
\end{align*}
Here, we used the computation
\begin{equation}\label{eqn:nabla-sigma-w-q-eta}
|\nabla^{\Sigma} (w^{-q} \eta)|^2 
%= |-q w^{-q-1} \eta \nabla^{\Sigma} w + w^{-q} \nabla^{\Sigma} \eta|^2 
= q^2 w^{-2(q+1)} \eta^2 |\nabla^{\Sigma} w|^2 + w^{-2q} |\nabla^{\Sigma} \eta|^2 - 2 q w^{-2q-1} \eta \la \nabla^{\Sigma} w, \nabla^{\Sigma} \eta \rg 
\end{equation}
Moreover, testing the weak inequality~\eqref{eqn:weak-supersolution-semilinear} with $\phi = w^{-2q-1} \eta^2$ implies
\begin{align*}
    - 2 \int_{\Sigma} w^{-(2q+1)} \eta \la \nabla^{\Sigma} \eta, \nabla^{\Sigma} w \rg + (2q+1) \int_{\Sigma} \eta^2 w^{-2(q+1)} |\nabla^{\Sigma} w|^2 \leq 4 \int_{\Sigma} w^{-2(q+1)} \eta^2. 
\end{align*}
We sum the two inequalities to obtain
\begin{align*}
    (2q+1-q^2) \int_{\Sigma} \eta^2 w^{-2(q+1)} |\nabla^{\Sigma} w|^2 \leq 2(1-q) \int_{\Sigma} w^{-2q-1} \eta \la \nabla^{\Sigma} w, \nabla^{\Sigma} \eta \rg + \int_{\Sigma} w^{-2q} |\nabla^{\Sigma} \eta|^2.
\end{align*}
For $q \in (0, 1 + \sqrt{2})$, we use $2q+1-q^2 > 0$ and Young's inequality to absorb the mixed term into
\[
2(1-q) \, w^{-2q-1} \eta \la \nabla^{\Sigma} w , \nabla^{\Sigma} \eta \rg \leq \tfrac{1}{2} (2q+1-q^2) \eta^2 w^{-2(q+1)} |\nabla^{\Sigma} w|^2 + C_q w^{-2q} |\nabla^{\Sigma} \eta|^2.
\]
Applying this result in the above integral estimate, we conclude that
\begin{equation}\label{eqn:int-sigma-w}
\int_{\Sigma} \eta^2  w^{-2(q+1)} |\nabla^{\Sigma} w|^2 \leq C_q \int_{\Sigma} w^{-2q} |\nabla^{\Sigma} \eta|^2.
\end{equation}
The same application of the stability inequality with $\phi = w^{-q} \eta$ as above, together with the bound
\[
|\nabla^{\Sigma} (w^{-q} \eta)|^2 \leq 2q^2 \eta^2 w^{-2(q+1)} |\nabla^{\Sigma} w|^2 + 2 w^{-2q} |\nabla^{\Sigma} \eta|^2
\]
due to the expression~\eqref{eqn:nabla-sigma-w-q-eta}, shows that the estimate~\eqref{eqn:int-sigma-w} implies
\begin{equation}\label{eqn:get-rid-of-it}
    \int_{\Sigma} w^{-2(q+1)} \eta^2 \leq C_q \int_{\Sigma} w^{-2q} |\nabla^{\Sigma} \eta|^2.
\end{equation}
Finally, we choose a non-negative cutoff function $\zeta$ satisfying $\zeta\equiv 1$ on $B^{\Sigma}_R(x_0)$, with $\textup{spt} \, \zeta \Subset B^{\Sigma}_{2R}(x_0)$ and $|\nabla^{\Sigma} \zeta| \leq C R^{-1}$.
We then apply the bound~\eqref{eqn:get-rid-of-it}, setting $p = 2(q+1)$ and $\eta = \zeta^{\frac{p}{2}}$, to find
\begin{align*}
    \int_{\Sigma} w^{-p} \zeta^{p} \leq C_p \int_{\Sigma} w^{-(p-2)} \zeta^{p-2} |\nabla^{\Sigma} \zeta|^2 \leq C_p \Bigl( \int_{\Sigma} w^{-p} \zeta^p \Bigr)^{\frac{p-2}{p}} \Bigl( \int_{\Sigma} |\nabla^{\Sigma} \zeta|^p \Bigr)^{\frac{2}{p}}
\end{align*}
where the second bound follows from H\"older's inequality.
Since $q \in (0,1+\sqrt{2})$, we have $p \in (2, 4+2 \sqrt{2})$.
The estimate~\eqref{eqn:L-p-estimate} for $w$ follows from rearranging terms and using the properties of $\zeta$.

Finally, when $\Sigma$ is a complete minimal hypersurface with $\textup{Vol}_{\Sigma} (B^{\Sigma}_{2R}(x_0)) \leq \Lambda R^{n-1}$, the conclusion is identical to the argument of~\cite{schoen-simon-yau}*{Theorem 3}.
Indeed, equation~\eqref{eqn:weak-supersolution-semilinear} implies that the function $u = 2 w^{-1}$ is a non-negative subsolution of
\begin{equation}\label{eqn:non-negative-subsolution}
    \Delta_{\Sigma} u + V u \geq 0, \qquad V = u^2.
\end{equation}
Moreover, for $n \leq 7$, the previous discussion shows that, for $p \in (n-1, 4 + 2 \sqrt{2})$, we have
\begin{equation}\label{eqn:u-lp-estimate}
\| u \|_{L^p(B^{\Sigma}_{3R/2})} \leq C(n,p,\Lambda) R^{\frac{n-1}{p}-1}.
\end{equation}
In particular, we obtain $V \in L^{\frac{p}{2}}$, for $\frac{p}{2} > \frac{n-1}{2}$, with a uniformly controlled scale-invariant norm, so the results of~\cite{morrey}*{Theorem 5.3.1} can be adapted to our situation, as explained in~\cite{schoen-simon-yau}*{Proof of Theorem 3, equations (2.11) - (2.12)}.
The Moser iteration described therein is available because $\Sigma$ has Euclidean volume growth and satisfies the Sobolev inequality of Hoffman-Spruck~\cite{hoffman-spruck}.

We combine the properties~\eqref{eqn:non-negative-subsolution} and~\eqref{eqn:u-lp-estimate}, so this Moser iteration argument shows that
\[
\esssup_ {B^{\Sigma}_R} u \leq C(n,\Lambda) R^{-1}, \qquad \text{where} \quad w := u^{-1}.
\]
Equivalently, we obtain $\essinf_{B^{\Sigma}_R} w \geq c(n,\Lambda)  R$.
This proves our second assertion.
\end{proof}

The dimensional threshold of Proposition~\ref{prop:stable-toda} is sharp.
Indeed, on $\mathbb{R}^{n-1}$, the radial function $w(x) = \frac{2}{\sqrt{n-2}} |x|$ is a distributional solution of~\eqref{eqn:weak-supersolution-semilinear} with $\frac{4}{w^2} = \frac{n-2}{|x|^2}$.
By the sharp Hardy inequality, the stability condition holds exactly when $n-2 \leq \tfrac{1}{4} (n-3)^2$; for integer $n \geq 1$, this amounts precisely to $n \geq 8$, exhibiting a singular solution of~\eqref{eqn:weak-supersolution-semilinear} that satisfies the stability inequality~\eqref{eqn:stability-semilinear}.
This phenomenon is intimately related to the stability of the axisymmetric cone in eight dimensions, which realizes the alternative of multiplicity-two convergence to a hyperplane as $s \uparrow 1$ in Theorem~\ref{thm:eight-dimensions}.

\section{Iteration for the sheet separation}\label{sec:iterate-sheet-separation}

We continue to work in the graphical setting of Section~\ref{sec:prelimiinaries-nonlocal-approximation}, decomposing the boundary $\partial E \cap \Omega$ of a stable $s$-minimal set as a union of graphs over $B^{n-1}_{1/2}$, given by functions $g_i$ satisfying~\eqref{eqn:slope-bound-g-i}.
We know from~\cite{nonlocal-approximation}*{Proposition~3.2} that
\begin{equation}\label{eqn:after-shrinking}
    (1-s)\int |g_{i+1} - g_i|^{-2} \, dx' \leq C_n
\end{equation}
for a dimensional constant $C_n$, after shrinking the cylinder by a fixed amount.

Our key observation in this Section is an iteration for the exponent in the above $L^{-2}$ separation estimate. 
This leads to the following result, which plays a pivotal role in our arguments:
\begin{proposition}\label{prop:improve-the-exponent}
There exist dimensional constants $\delta_*>0$, $s_*\in(0,1)$, $r_*>0$, and $C_n>0$ with the following property.
If $\delta \in (0,\delta_*)$ and $s \in (s_*,1)$, then every adjacent pair of graphs satisfies
\begin{equation}\label{eq:moment}
(1-s)^3 \int_{B_{r_*}^{n-1}} |g_{i+1} - g_i|^{-6} \, dx' \leq C_n
\end{equation}
provided that $|g_i(0)| \leq \frac{1}{4}$ and $|g_{i+1}(0)| \leq \frac{1}{4}$.
The constants are independent of the number of sheets and of the minimum separation between the sheets.
\end{proposition}

In view of Proposition~\ref{prop:stable-toda}, the exponent~$6$ in Proposition~\ref{prop:improve-the-exponent} could be improved to an $L^{-p}$ estimate for any $p<4+2\sqrt{2}$, upon taking $s$ sufficiently close to $1$ and $\delta$ sufficiently small.
However, this is refinement is not needed in the ensuing argument.

The proof of Proposition~\ref{prop:improve-the-exponent} uses an iteration argument through the nonlocal stability operator $\cL$ considered in Section~\ref{sec:prelimiinaries-nonlocal-approximation}, following ~\cite{nonlocal-approximation}*{Proposition 3.5} and the approach of~\cite{barrios-figalli-valdinoci}*{\S 3.1}.
We will prove this result at the end of this section.
The general strategy is similar to that of Proposition~\ref{prop:stable-toda}, applied to the function $w := \frac{g_{i+1} - g_i}{\sqrt{1-s}}$, but more terms appear in the nonlocal setting, and we must prove that their contribution is of lower order when $\delta \in (0,\delta_*)$ and $s \in (s_*,1)$.
Part of this control is Lemma~\ref{lemma:interaction-estimate}.
In addition, we obtain an integrated stability inequality with improved exponents, in the spirit of~\cite{ros-oton-serra}*{Lemmas 4.1 and 4.4} and~\cite{stable-thin-one-phase}*{Proposition 6.1}.
This step can be viewed as the nonlocal analogue of the Caccioppoli inequality derived from stability.

\begin{lemma}\label{lemma:xavi-and-joaquim}
Fix $q>0$ and $\ve>0$, and consider a bounded domain $D \subset \bR^{n-1}$ and a function $\eta \in C_c^1(D)$ with $\|\eta\|_{C^1} \leq M$.
Then, there exists a constant $C$, depending only on $D,n,q,\ve, C_0, M$, such that the following holds.
Let $K$ be a non-negative symmetric kernel on $D$, with
\[
K(x',y') \leq |x'-y'|^{-n-s}, \qquad \text{for} \; s \in (\tfrac{1}{2}, 1), \qquad x', y' \in D.
\]
Over $D$, let $\cL_K^D$ be the nonlocal operator associated to $K$ as in~\eqref{eqn:LK-operator-over-D}, and suppose that the smooth positive function $w$ satisfies
\begin{equation}\label{eqn:suppose-kernel-bound}
- (1-s) \, \cL_K^D w \leq \tfrac{41}{20} c_{n,s} \, w^{-s} + C_0 (1 + w).
\end{equation}
on the support of the function $\eta \in C_c^1(D)$.
Then, it holds that
\begin{equation}\label{eqn:integrating-stability}
    (1-s) \int_D (w^{-q} \eta) \, \cL_K^D (w^{-q} \eta) \leq \Bigl( \frac{41 q^2}{20(2q+1)} + \ve \Bigr) c_{n,s} \int_D w^{-2q-1-s} \eta^2 + C \int_D (w^{-2q} + \eta^2).
\end{equation}
Moreover, for functions $w_j, w_{j+1}$ and $\eta_j, \eta_{j+1}$ on $D$, we have the pointwise inequality
\begin{equation}\label{eqn:pointwise-inequality}
\begin{split}
    & \frac{w_j^{-2q} \eta_j^2 + w_{j+1}^{-2q} \eta_{j+1}^2}{(w_j + w_{j+1})^{1+s}} \\
    &\leq (1+\ve ) \, 2^{-1-s} (w_j^{-2q-1-s} \eta_j^2 + w_{j+1}^{-2q-1-s} \eta_{j+1}^2) + C \frac{w_j^{-2q} + w_{j+1}^{-2q}}{(w_j + w_{j+1})^{1+s}} |\eta_j - \eta_{j+1}|^2.
\end{split}
\end{equation}
\end{lemma}
\begin{proof}
The proof relies on two elementary inequalities, which we now establish.

\smallskip \noindent \textbf{Step 1.}
We observe that for any fixed $q>0$ and $\ve>0$ as above, there exists a $C = C(q,\ve)$ such that, for all $a,b>0$, all $\xi,v \in \bR$, and all $s \in (0,1)$, we have the inequalities
\begin{equation}\label{eqn:a-q-xi-b-q-v}
    (a^{-q} \xi - b^{-q} v)^2 + \bigl( \tfrac{q^2}{2q+1} + \ve \bigr) (a-b) ( a^{-2q-1} \xi^2 - b^{-2q-1} v^2) \leq C(a^{-2q}+b^{-2q})(\xi-v)^2
\end{equation}
as well as
\begin{equation}\label{eqn:ratio-(xi-v)2}
    \begin{split}
\frac{a^{-2q} \xi^2+b^{-2q} v^2}{(a+b)^{1+s}} &\leq (1+\ve) \, 2^{-1-s} (a^{-2q-1-s} \xi^2 + b^{-2q-1-s} v^2) + C \frac{a^{-2q} + b^{-2q}}{(a+b)^{1+s}} (\xi-v)^2.
    \end{split}
\end{equation}
Indeed, by symmetry we can assume that $a \geq b$.
For the first inequality, the fundamental theorem of calculus can be combined with the Cauchy-Schwarz inequality to show that
\begin{align*}
    (a^{-q} - b^{-q})^2 &= q^2 \Bigl( \int_b^a t^{-q-1} \, dt \Bigr)^2 \leq q^2 (a-b) \int_b^a t^{-2q-2} \,dt = \frac{q^2}{2q+1} (a-b)(b^{-2q-1} - a^{-2q-1}).
\end{align*}
Moreover, for $t \in [b,a]$, we have
\[
b^{-2q-1} - a^{-2q-1} = (2q+1) \int_b^a t^{-2q-2} \, dt \geq (2q+1) (a-b) a^{-2q-2}.
\]
Rearranging terms in this inequality gives
\[
(a-b)^2 a^{-2q-2} \leq \tfrac{1}{2q+1} (a-b) (b^{-2q-1} - a^{-2q-1}).
\]
Set $B := (a-b) (b^{-2q-1}- a^{-2q-1}) > 0$ for brevity.
Expanding $\xi = v + (\xi-v)$, we therefore obtain
\begin{align*}
    & (a^{-q} \xi - b^{-q} v)^2 + ( \tfrac{q^2}{2q+1} + \ve) \, (a-b) \, ( a^{-2q-1} \xi^2 - b^{-2q-1} v^2) \\
    & \leq - \ve B v^2 + C(q,\ve) \, \bigl( a^{-q} B^{\frac{1}{2}} |v| \, |\xi-v| + a^{-2q} (\xi-v)^2 \bigr) \\
    &\leq C(q,\ve) \, a^{-2q}(\xi-v)^2.
\end{align*}
In the last step, we used Young's inequality.
This proves~\eqref{eqn:a-q-xi-b-q-v}.

For the second bound, we use Chebyshev's and Jensen's inequalities, and convexity, to obtain
\[
\frac{a^{-2q-(1+s)} + b^{-2q-(1+s)}}{a^{-2q} + b^{-2q}} \geq \frac{a^{-(1+s)} + b^{-(1+s)}}{2} \geq \Bigl( \frac{a+b}{2} \Bigr)^{-(1+s)}.
\]
Notice that $1+s \in (1,2)$, due to $s \in (0,1)$.
Moreover, since $a \geq b$, we have
\[
\bigl| 1 - (1+\ve) ( \tfrac{a+b}{2a} )^{1+s} \bigr| \leq \max \{ 1, \ve \}.
\]
Thus, again expanding $\xi = v + (\xi-v)$, we arrive at
\begin{align*}
    & a^{-2q} \xi^2 + b^{-2q} v^2 - (1+\ve) ( \tfrac{a+b}{2} )^{1+s} (a^{-2q-1-s} \xi^2 + b^{-2q-1-s} v^2 ) \\
    &\leq - \ve\,  (a^{-2q} + b^{-2q}) \, v^2 + C(q,\ve) \, a^{-2q} \, ( |v| \, |\xi-v| + (\xi-v)^2) \\
    &\leq C(q,\ve) \, (a^{-2q} + b^{-2q}) (\xi-v)^2.
\end{align*}
In the last step, we again applied Young's inequality.
Dividing by $(a+b)^{1+s}$ proves the desired estimate~\eqref{eqn:ratio-(xi-v)2}, uniformly in $s \in (0,1)$.

\smallskip \noindent \textbf{Step 2.}
We now apply the inequalities of Step 1 as in~\cite{ros-oton-serra}*{Lemma 4.1 and (4.3)}.
Because $w$ is smooth and positive on the compact support of $\eta$, the evaluation of $\cL^D_{K} (w^{-q} \eta)$ is admissible.
Extending each cutoff function by zero $D$ and using the bound on the the $K$, we obtain
\begin{equation}\label{eqn:Kj-integral-bound}
    \begin{split}
    & (1-s) \sup_{x' \in D} \int_D |\eta(x') - \eta(y')|^2 \, K(x', y') \, dy' \\
    &\leq C(\eta)(1-s) \int_D |x'-y'|^2 \cdot |x' - y'|^{-n-s} \, dy' \\
    &\leq C(\eta) (1-s) \, |\bS^{n-2} |\int_0^{\textup{diam} \, D} r^{-s} \, dr \leq C' (\eta) \, ( \textup{diam} \, D)^{1-s},
    \end{split}
\end{equation}
uniformly in $j$ and in $s \in ( \frac{1}{2}, 1)$.
We apply the bound~\eqref{eqn:a-q-xi-b-q-v}, with $\ve$ replaced by $\frac{10}{41} \ve$, taking 
\[
(a,b,\xi,v) = ( w(x'), w(y'), \eta(x'), \eta(y')).
\]
Thus, we obtain the bound
\begin{equation}\label{eqn:wj-etaj-bound}
    \begin{split}
& ( w(x')^{-q} \eta(x') - w(y')^{-q} \eta(y'))^2 \\
&\leq - \bigl( \tfrac{q^2}{2q+1} + \tfrac{10}{41} \ve) ( w(x') - w(y')) [ w(x')^{-2q-1} \eta(x')^2 - w(y')^{-2q-1} \eta(y')^2] \\
& \quad + C (w(x')^{-2q} + w(y')^{-2q}) ( \eta(x') - \eta(y'))^2 \\
    \end{split}
\end{equation}
We integrate the inequality~\eqref{eqn:wj-etaj-bound} against the kernel $\frac{1-s}{2} K(x',y')$ and use symmetry in the expression~\eqref{eqn:Lu-kernel}, together with the preceding bound and the inequalities for $\cL^D_{K} w$ and for the kernel.

Using the expression~\eqref{eqn:Lu-kernel}, we compute the left-hand side under consideration as
\begin{align*}
    (1-s) \int_D (w^{-q} \eta) \, \cL^D_{K} (w^{-q} \eta) &= \frac{1-s}{2} \iint_{D \times D} (w(x')^{-q} \eta(x') - w(y')^{-q} \eta(y'))^2 K(x',y') \\
    &\leq - \Bigl( \frac{q^2}{2q+1} + \frac{10 \ve}{41} \Bigr) (1-s) \int_D w^{-2q-1} \eta^2 \cL^D_{K} w \\
    & \quad + C(1-s) \int_D w^{-2q}(x') \sup_{x' \in D} \int_D |\eta(x') - \eta(y')|^2 K(x',y') \, dy' \\
    &\leq  \Bigl( \frac{41 q^2}{20 (2q+1)} + \frac{\ve}{2} \Bigr) c_{n,s} \int_D w^{-2q-1-s} \eta^2 + C \int_D (w^{-2q} + w^{-2q-1} \eta^2).
\end{align*}
In the last step, we controlled the first integral on the right-hand side using the kernel bound~\eqref{eqn:suppose-kernel-bound}, and then estimated the second summand using~\eqref{eqn:Kj-integral-bound}.

Finally, we can apply Young's inequality to bound
\[
C(q,\ve,C_0) \, w^{-2q-1} \eta^2 \leq \tfrac{1}{2} c_{n,s} \ve \, w^{-2q-1-s} \eta^2 + C' \eta^2
\]
for a constant $C' = C'(q,\ve,C_0)$, which is uniform in $s \in ( \frac{1}{2},1)$.
Integrating this bound, we obtain
\begin{equation}\label{eqn:C-q-eps-young}
C(q,\ve,C_0) \int_D w^{-2q-1} \eta^2 \leq \frac{1}{2} c_{n,s} \ve \int_D w^{-2q-1-s} \eta^2 + C' \int_D \eta^2.
\end{equation}
Combining this term with the above integral estimate yields the inequality~\eqref{eqn:integrating-stability}.
Finally, the second bound~\eqref{eqn:pointwise-inequality} follows applying the bound~\eqref{eqn:ratio-(xi-v)2} of Step 1 with $(a,b,\xi,v) = (w_j, w_{j+1}, \eta_j, \eta_{j+1})$.
\end{proof}

We will prove Proposition~\ref{prop:first-prop} via an iteration argument, showing that an $L^{-2q}$ estimate in a larger ball implies an $L^{-2q-1-s}$ estimate in a smaller ball, when $q \leq \frac{41}{20}$.
This step can be viewed as an extension of~\cite{nonlocal-approximation}*{Proposition 3.2}, which proves such a property with $q=0$, and then deduces the $L^{-2}$ estimate by combining the $L^{-1-s}$ bound with a rough $L^{-\infty}$ control of the sheet separation, established in~\cite{nonlocal-approximation}*{Lemma 3.4}.

In our situation, the terms in the iteration toward $L^{-6}$ will involve the nonlocal interaction between the sheets.
We first show that this term can be absorbed, with a smaller coefficient, into the quantity being iterated.
The following argument proceeds similarly to~\cite{nonlocal-approximation}*{Lemma 3.3} and obtains an analogue of the bound (3.19) therein, which will now also involve $w_j^{-q}$.
A crucial modification is to control the integral by splitting the diagonal and off-diagonal contributions at a scale $\rho = L (g_{j+1} - g_j)(x')$, instead of the scale $\rho = [(g_{j+1} - g_j)(x')]^{\frac{1}{2}}$ used in~\cite{nonlocal-approximation}*{Lemma 3.3}.
Indeed, the natural expectation is that the two integrals become comparable at scale $\rho \asymp \sqrt{1-s}$, while Theorem~\ref{thm:sheet-happens} anticipates $\frac{g_{j+1} - g_j}{\sqrt{1-s}} \asymp c$, where $c$ may be very small.
This motivates the choice $\rho = L (g_{j+1} -g_j)(x')$ in the ensuing splitting argument, where $L = L_n$ may be as large as needed.

\begin{lemma}\label{lemma:interaction-estimate}
We can choose $\delta_* > 0$ such that for every $\delta \in (0,\delta_*)$, $s \in ( \frac{1}{2}, 1)$, $q \in [0,n]$, and $M < \infty$, there exists a constant $C = C(n,M)< \infty$ so that the following holds.
Suppose that the two sheets $\Gamma_j, \Gamma_{j+1}$ are graphs of functions $g_j < g_{j+1}$ over $B^{n-1}_R$, with $R \leq 2$ and
\[
\|D g_j \|_{L^{\infty}} + \|D^2 g_j \|_{L^{\infty}} < \delta, \qquad \|D g_{j+1} \|_{L^{\infty}} + \|D^2 g_{j+1} \|_{L^{\infty}} < \delta, \qquad |g_j(0)|, |g_{j+1}(0)| \leq \tfrac{1}{4}.
\]
Moreover, for some $r < \min \{ \frac{1}{4}, \frac{1}{2} R \}$, we assume that
\[
\int_{B^{n-1}_r} w^{-2q} \leq M, \qquad \text{where } \quad w := (1-s)^{- \frac{1}{1+s}} (g_{j+1} - g_j),
\]
and that $\eta$ is a non-negative cutoff function on $\bR^{n-1}$, compactly supported in $B^{n-1}_r$ and satisfying $0 \leq \eta \leq 1$ and $\|D \eta\|_{L^{\infty}} \leq M$.
Then, extending the function $w^{-q} \eta$ by zero outside $B^{n-1}_r$, we have
\begin{equation}\label{eqn:the-goal-interaction}
\begin{split}
    & (1-s) \iint_{\Gamma_j \times \Gamma_{j+1}} \frac{(w(x')^{-q} \eta(x') - w(y')^{-q} \eta(y'))^2}{|x-y|^{n+s}} \, d \cH^{n-1}_x \, d \cH^{n-1}_y \\
    &\leq \frac{c_{n,s}}{100} \int_{B^{n-1}_r} w^{-2q-1-s} \eta^2 + C_M (1+M).
\end{split}
\end{equation}
\end{lemma}
\begin{proof}
First, because $|D g_j| \leq \delta$, the graph Jacobians on $\Gamma_j, \Gamma_{j+1}$ are uniformly bounded by $1 + O(\delta)$.
Thus, we can transform the volume form $d \cH^{n-1}_x \, d \cH^{n-1}_y$ over $\Gamma_j \times \Gamma_{j+1}$ to $(1 + O(\delta)) \, dx' \, dy'$ over the $(n-1)$-dimensional base balls.
Moreover, for $(x,y) \in \Gamma_j \times (\Gamma_{j+1} \cap B_{\rho}^{n-1}(x'))$, we write $x = (x', g_j(x'))$ and $y = (y', g_{j+1}(y'))$ to obtain
\begin{align*}
    \bigl| ( g_{j+1}(y') - g_j(x')) - (g_{j+1} - g_j)(y') \bigr| &=  | g_j(y') - g_j(x')| \leq \delta \, |x'-y'| , \\
    \bigl| (g_{j+1}(y') - g_j(x')) - (g_{j+1} - g_j)(x') \bigr| &= |g_{j+1}(y') - g_{j+1}(x')| \leq \delta \, |x'-y'|.
\end{align*}
Thus, arguing as in Proposition~\ref{prop:first-prop} and~\cite{nonlocal-approximation}*{Lemma 3.3}, we can bound
\[
|x-y|^2 \geq \tfrac{1}{2} \max \{ |x'-y'|^2 + (g_{j+1}(x') - g_j(x'))^2 , |x'-y'|^2 + (g_{j+1}(y') - g_j(y'))^2 \}.
\]
Here, we used the simple fact that if $|\tilde{a} - a| \leq \delta t$, then
\begin{align*}
t^2 + \tilde{a}^2 = t^2 + a^2 + 2a( \tilde{a} - a) + ( \tilde{a} - a)^2 \geq t^2 + a^2 - 2 \delta  |a| t \geq (1-\delta) (t^2 + \tilde{a}^2)
\end{align*}
due to $2 \delta |a| t \leq \delta ( t^2 + a^2)$.
The above bound on $|x-y|^2$ follows from applying this inequality with 
\[
(t,\tilde{a}) = ( |x'-y'|, g_{j+1}(y') - g_j(x')), \qquad a \in \{ (g_{j+1} - g_j)(x'), (g_{j+1} - g_j)(y') \}.
\]
Hence, the interaction kernel $|x-y|^{-n-s}$ is bounded from above by $|x'-y'|^{-n-s}$, and by
\begin{equation}\label{eqn:interaction-kernel-bounds}
    |x-y|^{-n-s} \leq C_n \min \Bigl\{ \frac{1}{( |x'-y'|^2 + (g_{j+1}-g_j)(x')^2)^{\frac{n+s}{2}}}, \frac{1}{( |x'-y'|^2 + (g_{j+1}-g_j)(y')^2)^{\frac{n+s}{2}}} \Bigr\}.
\end{equation}
We will control the integral on the left-hand side of~\eqref{eqn:the-goal-interaction} as in~\cite{nonlocal-approximation}*{Lemma 3.3}.
We will repeatedly use the computation from the last step of that result, whereby, for all $s \in ( \frac{1}{2}, 1)$,
\begin{equation}\label{eqn:computation-last-step}
    \int_{\bR^{n-1} \setminus B^{n-1}_d (x') } |x-y|^{-(n+s)} \, d \cH^{n-1}_y \leq C_n |\bS^{n-2}| \int_d^{\infty} \bar{r}^{-2-s} \, d \bar{r} \leq C_n d^{-(1+s)}.
\end{equation}
This follows from applying the kernel bound $|x-y|^{-(n+s)} \leq |x'-y'|^{-(n+s)}$ and using polar coordinates on $y' \in \bR^{n-1}$, centered at $x'$, with volume form $\bar{r}^{n-2} \, d \bar{r} \, d \theta$ where $\bar{r} = |y'-x'|$.

Since $\eta = 0$ outside $B^{n-1}_r$, the contribution from pairs $(x',y') \in (B^{n-1}_r)^c \times (B^{n-1}_r)^c$ to the left-hand side of~\eqref{eqn:the-goal-interaction} is zero, so by symmetry, it suffices to restrict to $y' \in B^{n-1}_r$ in what follows.
Thus, by the Lipschitz bound $\| D \eta \|_{L^{\infty}} \leq M$ and the fact that $\eta = 0$ on $\partial B^{n-1}_r$, for $y' \in B^{n-1}_r$ we have
\begin{equation}\label{eqn:general-bound-eta}
    \begin{aligned}
        |\eta(x') - \eta(y')| &\leq \|D \eta \|_{L^{\infty}} |x'-y'|, \qquad & x' &\in B^{n-1}_r, \\
        |\eta(x') - \eta(y')| &\leq \|D \eta \|_{L^{\infty}} ( r - |y'|), \qquad & x' &\not\in B^{n-1}_r.
    \end{aligned}
\end{equation}
Notice that $\eta(x') = 0$ in the second case.
To estimate the kernel, we use the bound $|x-y| \geq |x'-y'|$ and observe that $(B^{n-1}_r)^c \subseteq B^{n-1}_{ r - |y'| }(y')$.
Thus, we can split the integration over $x'$ into
\[
B^{n-1}_r \cup (B^{n-1}_r)^c \subseteq B^{n-1}_R \cup ( \bR^{n-1} \setminus B^{n-1}_{r - |y'|}(y') ).
\]
Combining these properties, we therefore obtain
\begin{align*}
    & \iint_{\substack{\Gamma_j \times \Gamma_{j+1} \\ y' \in B^{n-1}_r}} \frac{(\eta(x') - \eta(y'))^2}{|x-y|^{n+s}} w(y')^{-2q} \, d \cH^{n-1}_x \, d \cH^{n-1}_y \\
    &\leq (1 + O(\delta)) \, C_n \, \|D \eta \|^2_{L^{\infty}} \int_{B^{n-1}_r} w(y')^{-2q} \left( \int_{B^{n-1}_R} |x'-y'|^{2-(n+s)} \, dx' \right) \, dy' \\
    & \quad + (1 + O(\delta)) \, C_n \, \|D \eta \|^2_{L^{\infty}} \int_{B^{n-1}_r} w(y')^{-2q} \cdot (r - |y'|)^2 \left( \int_{ \bR^{n-1} \setminus B^{n-1}_{r - |y'|}(y') } |x-y|^{-(n+s)} \, d \cH^{n-1}_x \right) \, dy'.
\end{align*}
For the integral in $x'$ in the first term, we use polar coordinates centered at $y'$ to estimate
\[
\int_{B^{n-1}_R} |x'-y'|^{2-(n+s)} \, dx' \leq C_n \, |\bS^{n-2}| \, (1-s)^{-1} R^{1-s} \leq C_n (1-s)^{-1}, \qquad \text{for } \; R \leq 2.
\]
For the second integral, in the complement of the ball $B^{n-1}_{r - |y'|}(y')$, we use the bound~\eqref{eqn:computation-last-step}, which produces a term $C_n ( r - |y'|)^{-(1+s)}$.
We then multiply by $(r - |y'|)^2$ and use the bound $r \leq 2$ to deduce that $C_n (r - |y'|)^{1-s} \leq C_n$.
Combining these estimates, we arrive at
\begin{equation}\label{eqn:joint-bound-first-term}
\begin{split}
    & (1-s) \iint_{  \substack{ \Gamma_j \times \Gamma_{j+1} \\ y' \in B^{n-1}_r } } \frac{ (\eta(x') - \eta(y'))^2}{|x-y|^{n+s}} w(y')^{-2q} \,  d \cH^{n-1}_x \, d \cH^{n-1}_y \leq C_nM^2 \int_{B^{n-1}_r} w^{-2q}.
\end{split}
\end{equation}
Thus, in what follows, we can consider points $(x,y) \in \Gamma_j \times \Gamma_{j+1}$ having $x',y' \in B^{n-1}_r$.

For some $L$ to be determined, and $\rho = L \, (g_{j+1} -g_j)(x')$, we will split the integral over $y \in \Gamma_{j+1}$ into two parts, one over the ball $B^{n-1}_{\rho}(x')$ and one over its complement.

\smallskip \noindent \textbf{Points $y'$ inside the ball $B^{n-1}_{\rho}(x')$.}
First, we add and subtract terms in the numerator to write
\begin{align*}
    w(x')^{-q} \eta(x') - w(y')^{-q} \eta(y') = \eta(x') \, [ w(x')^{-q} - w(y')^{-q}] + w(y')^{-q} \, ( \eta(x') - \eta(y')).
\end{align*}
Thus, we can use the inequality $(a+b)^2 \leq 2(a^2+b^2)$ to bound
\begin{equation}\label{eqn:(a+b)2-2a2b2}
    \begin{split}
    [ w(x')^{-q} \eta(x') - w(y')^{-q} \eta(y') ]^2 \leq 2 \, \eta(x')^2  \, (w(x')^{-q} - w(y')^{-q})^2 + 2 \, (\eta(x') - \eta(y'))^2 w(y')^{-2q}.
    \end{split}
\end{equation}
The integral of the second term is bounded by the estimate~\eqref{eqn:joint-bound-first-term}.
To control the first term of the above bound, we claim that, for $y \in \Gamma_{j+1} \cap B^{n-1}_{\rho}(x')$, it holds that
\begin{equation}\label{eqn:bound-in-the-ball-v2}
    (1-s) \int_{\Gamma_{j+1} \cap B_{\rho}^{n-1}(x')} \frac{(w(x')^{-q} -w(y')^{-q})^2 \eta(x')^2}{|x-y|^{n+s}} \, d \cH^{n-1}_y \leq C_n \delta^2 L^{\frac{1}{2}} w(x')^{-2q-1-s} \eta(x')^2.
\end{equation}
Indeed, for points $y' \in \Gamma_{j+1} \cap B^{n-1}_{\rho}(x')$, the gradient bounds~\eqref{eqn:gj-x-prime} and $|D(g_{j+1} - g_j)| \leq 2 \delta$ imply
\begin{align*}
    | (g_{j+1} - g_j)(x') -(g_{j+1} -g_j)(y') | &\leq 8 \delta^{\frac{1}{2}} \textstyle{\sqrt{(g_{j+1} - g_j)(x')}} \, |x'-y'| \leq 8 \delta^{\frac{1}{2}} \textstyle{\sqrt{(g_{j+1} - g_j)(x')}} \, \rho \\
    &\leq 8 \delta^{\frac{1}{2}} L \, [ (g_{j+1} - g_j)(x')]^{\frac{3}{2}}.
\end{align*}
Since $|g_j(0)|, |g_{j+1}(0)| \leq \frac{1}{4}$, the gradient bound on the $g_j$ also gives $|g_{j+1} - g_j| \leq C_n$.
Thus, after choosing $\delta>0$ sufficiently small, in terms of $L$, we arrive at
\begin{equation}\label{eqn:uniformly-comparable-bound}
    \tfrac{1}{2} (g_{j+1} - g_j)(x') \leq (g_{j+1} - g_j)(y') \leq \tfrac{3}{2} (g_{j+1} - g_j)(x'), \qquad \text{for every } \; y' \in B^{n-1}_{\rho}(x').
\end{equation}
In particular, this implies that the value of $(g_{j+1} - g_j)$ is uniformly comparable to $(g_{j+1}-g_j)(x')$ for every $y \in B_{\rho}^{n-1}(x')$.
Consequently, using the mean value theorem for the function $t \mapsto t^{-q}$, with $w = (1-s)^{-\frac{1}{1+s}} (g_{j+1} - g_j)$, and recalling that $0 \leq q \leq n$, we obtain
\begin{align*}
|w(x')^{-q} - w(y')^{-q}| &\leq \frac{2 q \delta \, (1-s)^{\frac{q}{1+s}} \, |x'-y'|}{ \min \{ (g_{j+1}-g_j)(x'), (g_{j+1}-g_j)(y') \}^{q+1} } \leq \frac{C_n \delta \, w(x')^{-q} \, |x'-y'|}{(g_{j+1} - g_j)(x')}.
\end{align*}
In the next step, we combine the kernel estimate~\eqref{eqn:interaction-kernel-bounds} with this bound to control the term~\eqref{eqn:bound-in-the-ball-v2} as
\begin{equation}\label{eqn:bound-in-the-ball}
\begin{split}
    & \int_{\Gamma_{j+1}  \cap B^{n-1}_{\rho}(x')} \frac{ (w(x')^{-q} - w(y')^{-q})^2 \eta(x')^2}{|x-y|^{n+s}} \, d \cH^{n-1}_y \\
    &\leq C_n \delta^2 \, w(x')^{-2q} \eta(x')^2 \int_{B^{n-1}_{\rho}(x')} \frac{|x'-y'|^2}{(g_{j+1}-g_j) (x')^2} \cdot \frac{dy'}{(|x'-y'|^2 + (g_{j+1}-g_j)(x')^2)^{\frac{n+s}{2}}}.
\end{split}
\end{equation}
Using polar coordinates $\bar{r}^{n-2} \, dr \, d \theta$ in $\mathbb{R}^{n-1}$ and changing variables to $t = \frac{\bar{r}}{|g_{j+1}(x') - g_j(x')|}$, we compute
\begin{align*}
& \int_{B^{n-1}_{\rho}(x')} \frac{|x'-y'|^2}{(g_{j+1}-g_j) (x')^2} \cdot \frac{dy'}{(|x'-y'|^2 + (g_{j+1}-g_j)(x')^2)^{\frac{n+s}{2}}} \\
&= \frac{|\bS^{n-2}|}{(g_{j+1} - g_j)(x')^{1+s}} \int_0^{\frac{\rho}{(g_{j+1}-g_j)(x')}} \frac{t^n}{(1+t^2)^{\frac{n+s}{2}}} \, dt = \frac{|\bS^{n-2}|}{(1-s) \, w(x')^{1+s}} \int_0^L t^n \, (1+t^2)^{-\frac{n+s}{2}} \, dt.
\end{align*}
The last integral is at most $3 \sqrt{L}$, for $s \in ( \frac{1}{2},1)$ and $L \geq L_n$ large.
Indeed, we can split the integral into $t \in [0,1]$ and $t \geq 1$, as well as bound $t^n (1+t^2)^{- \frac{n+s}{2}} \leq t^{-s} \leq t^{- \frac{1}{2}}$ for $t \geq 1$, to estimate
\begin{align*}
    \int_0^L t^n (1+t^2)^{- \frac{n+s}{2}} \, dt &\leq \int_0^1 t^n \, dt + \int_1^L t^{- \frac{1}{2}} \leq C_n + 2 L^{\frac{1}{2}} \leq 3 \sqrt{L}.
\end{align*}
The asserted inequality~\eqref{eqn:bound-in-the-ball-v2} follows from applying this bound to~\eqref{eqn:bound-in-the-ball}.

Finally, we apply the inequality~\eqref{eqn:(a+b)2-2a2b2} and integrate the bound~\eqref{eqn:bound-in-the-ball} over the set of $x \in \Gamma_j \cap (B^{n-1}_r \times \bR)$, using also the estimate~\eqref{eqn:joint-bound-first-term}.
Hence, we arrive at
\begin{equation}\label{eqn:final-in-the-ball}
    \begin{split}
    & (1-s) \iint_{ \substack{ \Gamma_j \times ( \Gamma_{j+1} \cap B^{n-1}_{\rho}(x')) \\ x',y' \in B^{n-1}_r } } \frac{ ( w(x')^{-q} \eta(x') - w(y')^{-q} \eta(y'))^2 }{|x-y|^{n+s}} \, d \cH^{n-1}_x \, d \cH^{n-1}_y \\
    &\leq C_n \delta^2 L^{\frac{1}{2}} \int_{B^{n-1}_r} w^{-2q-1-s} \eta^2 + C_n M^2 \int_{B^{n-1}_r} w^{-2q}.
    \end{split}
\end{equation}

\noindent \textbf{Points $y'$ outside the ball $B^{n-1}_{\rho}(x')$.}
For $y \in \Gamma_{j+1} \setminus B^{n-1}_{\rho}(x')$, we bound terms using
\begin{equation}\label{eqn:difference-bound-crude}
    (w(x')^{-q} \eta(x') - w(y')^{-q} \eta(y'))^2 \leq 2 w(x')^{-2q} \eta(x')^2 + 2 w(y')^{-2q} \eta(y')^2
\end{equation}
and $|x-y|^{-(n+s)} \leq |x'-y'|^{-(n+s)}$.
Moreover, combining the inequalities $|D(g_{j+1} - g_j)| \leq 2 \delta$ and $|y'-x'| > \rho = L (g_{j+1} - g_j)(x')$ in this region, and imposing $\delta < (10 L)^{-1}$, we obtain
\[
(g_{j+1} - g_j)(y') \leq (g_{j+1} - g_j)(x') + 2 \delta \, |x'-y'| < (L^{-1} + 2 \delta) \, |x'-y'| \leq \tfrac{3}{2} L^{-1} |x'-y'|.
\]
Thus, $|y'-x'| > L \cdot \max \{ (g_{j+1}-g_j)(x'), \frac{2}{3} (g_{j+1} - g_j)(y') \}$.
We will also use the fact that 
\[
\rho^{1+s} = L^{1+s} ( g_{j+1} - g_j)(x')^{1+s} = L^{1+s} \, (1-s) \, w(x')^{1+s}.
\]
We integrate the bound~\eqref{eqn:difference-bound-crude} first in the variable absent from the numerator, using~\eqref{eqn:computation-last-step} to bound
\begin{align*}
    (1-s)\int_{\Gamma_{j+1} \setminus B^{n-1}_{\rho}(x')} |x-y|^{-(n+s)} \, d \cH^{n-1}_y &\leq C_n (1-s)\rho^{-(1+s)} = C_n L^{-1-s} w(x')^{-1-s}.
\end{align*}
Thus, writing $L^{-1-s} \leq L^{-\frac{3}{2}}$ for $s \in ( \frac{1}{2}, 1)$, we obtain the integral estimate
\begin{align*}
    & (1-s) \iint_{\Gamma_j \times ( \Gamma_{j+1} \setminus B^{n-1}_{\rho}(x'))} \frac{w(x')^{-2q} \eta(x')^2}{|x-y|^{n+s}} \, d \cH^{n-1}_x \, d \cH^{n-1}_y \\ 
    &\leq C_n L^{-\frac{3}{2}} \int_{B^{n-1}_r} w(x')^{-2q} \eta(x')^2 w(x')^{-1-s} \, d x' = C_n L^{-\frac{3}{2}} \int_{B^{n-1}_r} w^{-2q-1-s} \eta^2 .
\end{align*}
Next, for the integral involving $w(y')^{-2q} \eta(y')^2$, we use the previous observation that $|x'-y'| > \frac{2}{3} L (g_{j+1} - g_j)(y')$ whenever $y' \not\in B^{n-1}_{\rho}(x')$.
Integrating first in $x'$, we therefore obtain
\begin{align*}
    & (1-s) \iint_{\substack{\Gamma_j \times ( \Gamma_{j+1} \setminus B^{n-1}_{\rho}(x')) \\ x',y' \in B^{n-1}_r}} \frac{w(y')^{-2q} \eta(y')^2}{|x-y|^{n+s}} \, d \cH^{n-1}_x \, d \cH^{n-1}_y \\
    &\leq C_n (1-s) \int_{B^{n-1}_r} w(y')^{-2q} \eta(y')^2 \left( \int_{ \{ |x'-y'| > \frac{1}{2} L (g_{j+1} - g_j)(y') \} } |x'-y'|^{-(n+s)} \, dx'  \right) \, dy' \\
    &\leq C_n L^{- \frac{3}{2}} \int_{B_r^{n-1}} w(y')^{-2q} \eta(y')^2 w(y')^{-1-s} \, dy' = C_n L^{- \frac{3}{2}} \int_{B^{n-1}_r} w^{-2q-1-s} \eta^2.
\end{align*}
Here, we used the bound~\eqref{eqn:computation-last-step} to estimate the integral in the second step by $C_n L^{-1-s} (g_{j+1} - g_j)(y')^{-1-s}$, and hence by $\leq C_n L^{- \frac{3}{2}} (1-s)^{-1} w(y')^{-1-s}$, for $s \in ( \frac{1}{2}, 1)$.
We therefore obtain
\begin{equation}\label{eqn:combine-these-bounds-exterior}
\begin{split}
    & (1-s) \iint_{ \substack{ \Gamma_j \times ( \Gamma_{j+1} \setminus B^{n-1}_{\rho}(x')) \\ x',y' \in B^{n-1}_r } } \frac{(w(x')^{-q} \eta(x') - w(y')^{-q} \eta(y'))^2}{|x-y|^{n+s}} \, d \cH^{n-1}_x \, d \cH^{n-1}_y  \\
    &\leq C_n L^{- \frac{3}{2}} \int_{B^{n-1}_r} w^{-2q-1-s} \eta^2.
\end{split}
\end{equation}

\noindent \textbf{Conclusion of the proof of Lemma~\ref{lemma:interaction-estimate}}.
We now combine the bounds~\eqref{eqn:joint-bound-first-term}, ~\eqref{eqn:final-in-the-ball}, and~\eqref{eqn:combine-these-bounds-exterior}.
Moreover, for every $\ve>0$, Young's inequality produces a constant $C(\ve)< \infty$ such that
\begin{equation}\label{eqn:young-inequality-bound}
    w^{-2q-1} \leq \ve \, w^{-2q-1-s} + C(\ve) \, w^{-2q}.
\end{equation}
The constant $C(\ve)$ is uniform in $s \in ( \frac{1}{2},1)$.
We therefore obtain
\begin{equation}\label{eqn:desired-bound}
    \begin{split}
    & (1-s) \iint_{\Gamma_j \times \Gamma_{j+1}} \frac{(w(x')^{-q} \eta(x') - w(y')^{-q} \eta(y'))^2}{|x-y|^{n+s}} \, d \cH^{n-1}_x \, d \cH^{n-1}_y \\
    &\leq (C_n \delta^2 L^{\frac{1}{2}} + C_n L^{- \frac{3}{2}}) \int_{B^{n-1}_r} w^{-2q-1-s} \eta^2 + C_n M^2 \int_{B^{n-1}_r} w^{-2q}.
    \end{split}
\end{equation}
The quantity $c_{n,s}$ appearing in~\cite{nonlocal-approximation} is $c_{n,s} := |\bS^{n-2}| \int_0^{\infty} t^{n-2}(1+t^2)^{- \frac{n+s}{2}} \, dt$, and is bounded below by a positive dimensional constant $c_* := \inf_{s \in ( \frac{1}{2}, 1)} c_{n,s}>0$ for $s \in ( \frac{1}{2},1)$.
Thus, we can first choose $L = L(n) \geq 1$ sufficiently large so that $C_n L^{- \frac{3}{2}} \leq 10^{-3} c_*$, and then choose $\delta_* = \delta_*(n)>0$ sufficiently small so that $\delta_* < (10 L)^{-1}$ and $C_n \delta_*^2 L^{\frac{1}{2}} < 10^{-3} c_*$.
Applying these properties in the bound~\eqref{eqn:desired-bound} establishes the desired inequality~\eqref{eqn:the-goal-interaction}, proving our assertion.
\end{proof}

We now combine these estimates to prove Proposition~\ref{prop:first-prop}.
We obtain the desired estimate through an iteration that shares some similarities with the one implemented by Fazly~\cite{fazly}*{\S 4} in the setting of stable solutions to nonlocal elliptic systems.

\begin{proof}[Proof of Proposition~\ref{prop:improve-the-exponent}]
We will perform a discrete iteration in three steps.
For abbreviation, we introduce the normalized function $w_i := (1-s)^{-\frac{1}{1+s}} (g_{i+1} - g_i)$, where
\[
\tfrac{1}{2} \textstyle{\sqrt{1-s}} < (1-s)^{\frac{1}{1+s}} < \textstyle{\sqrt{1-s}} \qquad \text{for } \; s \in (0,1).
\]
We fix $0<b<B<\frac{1}{2}$ and a sufficiently small ball $B^{n-1}_r$.
We will prove that a uniform $L^{-2q}$ bound
\begin{equation}\label{eqn:inductive-hypothesis}
\int_{B^{n-1}_r} w_j^{-2q} \leq M, \qquad \text{whenever} \quad |g_j(0)|, |g_{j+1}(0)| \leq B,
\end{equation}
implies the improved $L^{-2q-1-s}$ bound, in a smaller ball,
\begin{equation}\label{eqn:iterated-bound}
    \int_{B^{n-1}_{r/4}} w_i^{-2q-1-s} \leq C(1+M), \qquad \text{whenever} \quad |g_i(0)|, |g_{i+1}(0)| \leq b
\end{equation}
for the successive choices of $q \in \{ 1, \frac{3}{2}, \frac{41}{20} \}$.
The constants here may depend on the fixed radii and on the heights $b,B$, but not on the number $N$ of sheets or on the separation between them.

We work on the larger ball $B^{n-1}_r$, using the inductive hypothesis~\eqref{eqn:inductive-hypothesis}.
We first localize the kernel $K_j$ to this ball: to distinguish the localized kernel, let us write, as in Lemma~\ref{lemma:xavi-and-joaquim},
\[
\cL^r_{K_j} u(x') := \textup{PV} \int_{B^{n-1}_r} ( u(x') - u(y')) \, K_j(x', y') \, dy'.
\]
This definition corresponds to $D = B^{n-1}_r$.
The restriction of $K_j$ to $B^{n-1}_r$ retains the symmetry and kernel bounds from Proposition~\ref{prop:first-prop}.
To estimate the far-field contribution, for $x' \in B^{n-1}_{r/2}$ we let
\[
a_{j,r}(x') := \int_{B^{n-1}_{1/2} \setminus B^{n-1}_{r/2}} K_j(x',y') \, dy'
\]
Since $|x'-y'| \geq \frac{r}{2}$ in this integral, and $K_j(x',y') \leq |x'-y'|^{-(n+s)}$, the estimate~\eqref{eqn:computation-last-step} implies that $a_{j,r}(x') \leq C_n r^{-(1+s)}$.
Moreover, the localized operator $\cL^r_{K_j}$ is related to the original operator by
\[
\cL^r_{K_j} = \cL_{K_j} w_j(x') - \int_{B^{n-1}_{1/2} \setminus B_r^{n-1}} (w_j(x') - w_j(y')) \, K_j(x',y') \, dy'.
\]
Combining this property with the bound~\eqref{eqn:suppose-kernel-bound} and using $w_j = (1-s)^{- \frac{1}{1+s}} (g_{j+1} - g_j)$, we arrive at
\[
- (1-s) \, \cL^r_{K_j} w_j \leq \tfrac{41}{20} c_{n,s} \, w_j^{-s} + C_0 + (1-s) C_n r^{-1-s} w_j 
\]
In the ensuing iteration, we can ensure that $r \geq 10^{-10(n+1)}$ for each resulting radius.
Hence, $(1-s)r^{-1-s} \leq (1-s) 10^{20(n+1)} \leq 1+C_0$ for $s \in (1 - \frac{10^{-20(n+1)}}{1+C_0}, 1)$.
Thus, the above bounds become
\[
K_j(x',y') \leq |x'-y'|^{-(n+s)}\]and\[ - (1-s) \, \cL^r_{K_j} w_j \leq \tfrac{41}{20} c_{n,s} w_j^{-s} + (1+C_0) (1+ w_j).
\]
In particular, the conditions of Lemma~\ref{lemma:xavi-and-joaquim} hold in our situation, with $D = B^{n-1}_r$.
To combine this result with the stability inequality of Proposition~\ref{prop:first-prop}, we first notice that
\begin{equation*}
    \begin{split}&
\frac{1-s}{(g_{j+1} - g_j)^{1+s}} = w_j^{-1-s}, \qquad \frac{1-s}{(g_{j+1} -  g_{j-1})^{1+s}} = (w_{j-1} + w_j)^{-1-s}, \\&
\qquad\qquad{\mbox{and}}\qquad \frac{1-s}{(g_{j+2} - g_j)^{1+s}} = ( w_j + w_{j+1})^{-1-s}.
\end{split}\end{equation*}
As in Proposition~\ref{prop:first-prop}, we denote by $\Gamma_j$ the graphical sheets of $\partial E \cap \Omega$, namely
\[
\Gamma_j := \{ (x', g_j(x')) : x' \in B^{n-1}_{1/2}, \; -1 < g_j(x') < 1 \}.
\]
We apply~\eqref{eqn:stability-iteration} with $\phi = w_j^{-q} \eta_j$, multiply both sides by $(1-s)$, and use~\eqref{eqn:integrating-stability} to obtain
\allowdisplaybreaks{
\begin{equation}\label{eqn:iteration-first-step}
\begin{split}
& \left[ 2 - C \delta - (1+C \delta^2) \Bigl( \frac{41 q^2}{10(2q+1)} + \ve \Bigr) \right] \, c_{n,s} \int_{B^{n-1}_r} w_j^{-2q-1-s} \eta_j^2 \, dx ' \\
&\leq \frac{1-s}{2} \iint_{\Gamma_j \times \Gamma_{j+1}} \frac{(w_j(x')^{-q} \eta_j(x') - w_j(y')^{-q} \eta_j(y'))^2}{|x-y|^{n+s}} \, d \cH^{n-1}_x \, d \cH^{n-1}_y \\
& \quad + \Bigl( \frac{1}{2} + C \delta \Bigr) \, c_{n,s} \int_{B^{n-1}_r} \Bigl[ \frac{1}{(w_{j-1} + w_j)^{1+s}} + \frac{1}{(w_j + w_{j+1})^{1+s}} \Bigr] \, w_j^{-2q} \, \eta_j^2 \, dx' \\
& \quad + C(n,q,\ve,\eta_j) \int_{B^{n-1}_r} (w^{-2q}_j + \eta_j^2) \, dx'.
\end{split}
\end{equation}}
The remainder $C(1-s+\delta^2) \int_{B^{n-1}_r} w_j^{-2q} \eta_j^2$ from~\eqref{eqn:stability-iteration} was absorbed into the last integral using the bound $\|\eta_j\|_{C^1} \leq M$.
The final constant $C$ depends on $n,q,\ve,M$, but not on $j,N,s$ or the minimum gap.
The last summand in~\eqref{eqn:integrating-stability} only occurs for indices $j$ with $\eta_j \not\equiv 0$, else the term in~\eqref{eqn:C-q-eps-young} is zero.

In the next step, we will localize the estimate in a distinguished sheet index $i$, as in~\cite{nonlocal-approximation}*{Proposition 3.2}.
For some sufficiently small $\kappa>0$, to be determined, we multiply the inequality~\eqref{eqn:iteration-first-step} for $j$ by $e^{- \kappa |j-i|}$ and sum over all sheets.
Notice that, independently of the number $N$,
\[
\sum_{j=1}^{N-1} e^{- \kappa |j-i|} \leq 1 + \frac{2}{e^{\kappa}-1} =: C_{\kappa}, \qquad e^{- \kappa} \leq \frac{e^{- \kappa |(j+1)-i|}}{e^{- \kappa |j-i|}} \leq e^{\kappa}.
\]
Summing the second terms of the right-hand side, which consist of neighboring potentials, gives
\allowdisplaybreaks{
\begin{equation}\label{eqn:cross-term-to-estimate}
    \begin{split}
    & \sum_{j=1}^{N-1} e^{- \kappa |j-i|} \int_{B^{n-1}_r} \Bigl[ \frac{1}{(w_{j-1} + w_j)^{1+s}} + \frac{1}{(w_j + w_{j+1})^{1+s}} \Bigr] \, w_j^{-2q} \eta_j^2 \\
    &= \sum_{j=1}^{N-2} \int_{B^{n-1}_r} \frac{e^{- \kappa |j-i|} w_j^{-2q} \eta_j^2 + e^{- \kappa |(j+1)-i|} w_{j+1}^{-2q} \eta_{j+1}^2}{(w_j + w_{j+1})^{1+s}}.
    \end{split}
\end{equation}}
For each summand, we bound both exponential weights using
\[
\max \{ e^{- \kappa |j-i|} , e^{- \kappa |j+1 - i|} \} \leq e^{\kappa} \, \min \{ e^{- \kappa |j-i|} , e^{- \kappa |j+1-i|} \}.
\]
Thus, we can apply the bound~\eqref{eqn:pointwise-inequality}, allowing for an additional term $e^{\kappa}$, to obtain
\begin{equation}\label{eqn:omega-j-bound}
\begin{split}
& \frac{e^{- \kappa |j-i|} w_j^{-2q} \eta_j^2 + e^{- \kappa |j+1-i|} w_{j+1}^{-2q} \eta_{j+1}^2}{(w_j + w_{j+1})^{1+s}} \\
&\leq (1+\ve) \, e^{\kappa} \, 2^{-1-s} ( e^{- \kappa |j-i|} w_j^{-2q-1-s} \eta_j^2 + e^{- \kappa |j+1-i|} w_{j+1}^{-2q-1-s} \eta_{j+1}^2) \\
& \quad + C ( e^{- \kappa |j-i|} + e^{- \kappa |j+1-i|} ) \frac{w_j^{-2q} + w_{j+1}^{-2q}}{( w_j + w_{j+1})^{1+s}} |\eta_j - \eta_{j+1}|^2.
\end{split}
\end{equation}
We specialize the computation with test functions $\eta_j$ chosen as in~\cite{nonlocal-approximation}*{Proposition 3.2},
\[
\eta_j(x') := \eta(x') \zeta( g_j(x')) \zeta (g_{j+1}(x')).
\]
Here, $0 \leq \eta, \zeta \leq 1$ are cutoff functions satisfying
\begin{align*}
\eta = 1 \quad \text{on } \; B^{n-1}_{r/4}, \qquad \text{spt} \, \eta \Subset B^{n-1}_{r/2}, \qquad \textup{spt} \, \zeta \Subset ( - \tfrac{b+B}{2}, \tfrac{b+B}{2}) ,
\end{align*}
with $\zeta=1$ on a neighborhood of $[-b,b]$.

For $\delta$ sufficiently small, this construction ensures that $\eta_i = 1$ on the ball $B^{n-1}_{r/4}$, where $i$ is the distinguished index with $|g_i(0)|, |g_{i+1}(0)| \leq b$.
Moreover, our construction ensures that if $\eta_j \not \equiv 0$, then there is some $x' \in B^{n-1}_{r/2}$ at which $\zeta(g_j(x')) \zeta(g_{j+1}(x')) \neq 0$, and hence
\[
|g_j(x')| < \tfrac{1}{2} (b+B), \qquad |g_{j+1}(x')| < \tfrac{1}{2} (B+b).
\]
Using the slope bound for the functions $g_{\ell}$, with $\ell \in \{ j, j+1\}$, we therefore obtain
\[
|g_{\ell}(0)| \leq |g_{\ell}(x')| + \delta \, |x'| < \tfrac{1}{2} (b+B) + \tfrac{1}{2} \delta r < B
\]
upon choosing $\delta r < \frac{B-b}{10}$.
Thus, we have $|g_j(0)|, |g_{j+1}(0)|\leq B$, so our inductive hypothesis~\eqref{eqn:inductive-hypothesis} applies.
In particular, we deduce that
\begin{equation}\label{eqn:bounded-whenever}
\int_D w^{-2q}_j \leq M \qquad \text{whenever} \quad \eta_j \not\equiv 0.
\end{equation}
Consequently, for $r<1$, we can sum these terms and use $\eta_j \leq 1$ to obtain
\begin{equation}\label{eqn:sum-exponential-terms}
    \sum_{\substack{ 1 \leq j \leq N-1 \\ \eta_j \not \equiv 0 }} e^{- \kappa |j-i|} \int_{B^{n-1}_r}(w_j^{-2q} + \eta_j^2) \leq ( M + |B^{n-1}_r| ) \sum_{j=1}^{N-1} e^{- \kappa |j-i|} \leq C(\kappa, B)(1+M)
\end{equation}
Moreover, working as in Step 2 of~\cite{nonlocal-approximation}*{Proposition 3.2}, we observe that
\begin{align*}
& \eta_j(x') - \eta_{j+1}(x') = \eta(x') \,  \zeta(g_{j+1}) \, ( \zeta(g_j ) - \zeta(g_{j+2}) ), \\
& \implies |\eta_j - \eta_{j+1}| \leq C \, \min \{ g_{j+2} - g_j, 1 \}.
\end{align*}
Combining this property with the fact that~$(w_j + w_{j+1})^{1+s} = \frac{(g_{j+2} - g_j)^{1+s}}{1-s}$, we deduce that
\begin{equation}\label{eqn:bound-eta-j}
    \frac{|\eta_j - \eta_{j+1}|^2}{(w_j + w_{j+1})^{1+s}} \leq C(1-s) \frac{\min \{ (g_{j+2} - g_j)^2, 1 \} }{ (g_{j+2} - g_j)^{1+s}} \leq C(1-s).
\end{equation}
Here, we used the fact that $\min \{ t^2, 1 \} \, t^{-1-s} \leq 1$ for $t>0$ and $s \in (0,1)$.

As above, we can obtain a further bound when $\eta_j, \eta_{j+1}$ are not both identically zero.
Indeed, if $\eta_j \not\equiv 0$, then out previous computation shows that
\[
|g_j(0)| < \tfrac{1}{2} (b+B) + \tfrac{1}{20} (B-b), \qquad |g_{j+1}(0)| < \tfrac{1}{2} (b+B) + \tfrac{1}{20} (B-b)
\]
provided that $\delta r< \frac{B-b}{10}$.
If $|g_{j+2}(0)| > B$, then 
\[
g_{j+2}(0) - g_{j+1}(0) > B - \tfrac{1}{2} (b+B) - \tfrac{1}{20} (B-b) > \tfrac{1}{5} (B-b).
\]
Hence, we may again use the slope bounds for $g_{j+1}, g_{j+2}$ to deduce that
\[
(g_{j+2} - g_{j+1})(x') > \tfrac{9}{20} (B-b) - 2\delta \, |x'| > \tfrac{9}{20}(B-b) - \tfrac{1}{5}(B-b) = \tfrac{1}{4}(B-b)
\]
for every $x' \in B^{n-1}_r$.
We therefore obtain
\begin{equation}\label{eqn:w-j+1}
w_{j+1}^{-2q} = (1-s)^{\frac{2q}{1+s}} (g_{j+2} - g_{j+1})^{-2q} \leq C(b,B).
\end{equation}
The argument bounding $w_j^{-2q}$ in this case is analogous.
Consequently, using the bounds~\eqref{eqn:bounded-whenever}, \eqref{eqn:bound-eta-j}, and~\eqref{eqn:w-j+1}, we conclude that
\begin{equation}\label{eqn:sum-separate-terms}
\begin{split}
& \sum_{j=1}^{N-2} (e^{- \kappa |j-i|} + e^{- \kappa |j+1-i|}) \int_{B^{n-1}_r} \frac{w_j^{-2q} + w_{j+1}^{-2q}}{(w_j + w_{j+1})^{1+s}} \, |\eta_j - \eta_{j+1}|^2 \\
&\leq C_{\kappa}(1-s) \sum_{\substack{1 \leq j \leq N-2 \\ \eta_j, \eta_{j+1} \; \text{not both $\equiv 0$}}} ( e^{- \kappa |j-i|} + e^{- \kappa |j+1-i|}) \int_{B^{n-1}_r} (w^{-2q}_j + w^{-2q}_{j+1}) \\
&\leq C(\kappa,b,B)(1-s)(1+M) \leq C(\kappa,b,B) (1+M).
\end{split}
\end{equation}
Moreover, each term in the expression~\eqref{eqn:omega-j-bound} occurs in at most two adjacent pairs.
Thus, we can sum the inequality~\eqref{eqn:cross-term-to-estimate} and apply it to~\eqref{eqn:iteration-first-step}, followed by using the bounds~\eqref{eqn:sum-exponential-terms} and~\eqref{eqn:sum-separate-terms}.
We obtain
\begin{equation}\label{eqn:sum-adjacent-reciprocals}
    \begin{split}
& \sum_{j=1}^{N-1} e^{- \kappa |j-i|} \int_{B^{n-1}_r} \Bigl[ \frac{1}{(w_{j-1} + w_j)^{1+s}} + \frac{1}{(w_j + w_{j+1})^{1+s}} \Bigr] \, w_j^{-2q} \eta_j^2 \\
& \leq (1+\ve) e^{\kappa} 2^{-s} \sum_{j=1}^{N-1} e^{- \kappa |j-i|} \int_{B^{n-1}_r} w_j^{-2q-1-s} \eta_j^2 \\
& \quad + C \sum_{j=1}^{N-2} ( e^{- \kappa |j-i|} + e^{- \kappa |j+1-i|}) \int_{B^{n-1}_r} \frac{w_j^{-2q} + w_{j+1}^{-2q}}{(w_j + w_{j+1})^{1+s}} |\eta_j - \eta_{j+1}|^2 \\
&\leq (1+\ve) e^{\kappa} 2^{-s} \sum_{j=1}^{N-1} e^{- \kappa |j-i|} \int_{B^{n-1}_r} w_j^{-2q-1-s} \eta_j^2 + C(\kappa,b,B) (1+M).
    \end{split}
\end{equation}
Next, we insert the inequality~\eqref{eqn:sum-adjacent-reciprocals} into the previous bound~\eqref{eqn:iteration-first-step} and rearrange the first term of the right-hand side to the left.
This yields
\begin{align*}
& \left[ 2 - C \delta - (1+C \delta^2) \Bigl( \frac{41 q^2}{20 (2q+1)} + \ve \Bigr) - \Bigl( \frac{1}{2} + C \delta \Bigr) (1+\ve) e^{\kappa} 2^{-s} \right] \, c_{n,s} \sum_{j=1}^{N-1} e^{- \kappa |j-i|} \int w_j^{-2q-1-s} \eta_j^2 \\
&\leq \frac{1-s}{2} \sum_{j=1}^{N-1} e^{- \kappa |j-i|} \iint_{\Gamma_j \times \Gamma_{j+1}} \frac{ ( w_j(x')^{-q} \eta_j(x') - w_j(y')^{-q} \eta_j(y'))^2 }{|x-y|^{n+s}} \, d \cH^{n-1}_x \, d \cH^{n-1}_y + C(1+M).
\end{align*}
For $q \in \{ 1, \frac{3}{2}, \frac{41}{20} \}$, we can choose $\ve, \kappa, \delta_*>0$ sufficiently small and $s_* \in (0,1)$ sufficiently close to $1$ so that, for $\delta \in (0,\delta_*)$ and $s \in (s_*,1)$, the coefficient of the left-hand side becomes 
\begin{equation}\label{eqn:q-delta-epsilon}
2 - C \delta - (1+C \delta^2) \bigl( \tfrac{41}{20} \tfrac{q^2}{2q+1} + \ve \bigr) - ( \tfrac{1}{2} + C \delta ) (1+\ve) e^{\kappa} 2^{-s} \geq (2 - \tfrac{1}{4}) - \tfrac{1}{100} - \tfrac{41}{20} \tfrac{q^2}{2q+1}.  
\end{equation}
For $q \leq \frac{41}{20}$, we compute that
\[
\Bigl( 2 - \frac{1}{4} \Bigr) - \frac{1}{100} - \frac{41}{20} \frac{q^2}{2q+1} \geq 2- \frac{1}{4}- \frac{1}{100} - \frac{41}{20} \frac{( \frac{41}{20})^2}{2 \cdot \frac{41}{20}+1} = \frac{2071}{40800} > \frac{1}{20}.
\]
Moreover, the conditions of Lemma~\ref{lemma:interaction-estimate} are satisfied for the function $w_j$ and the interaction of the sheets $(\Gamma_j , \Gamma_{j+1})$.
Consequently, the bound~\eqref{eqn:the-goal-interaction} shows that the above right-hand side is bounded by
\begin{align*}
& \frac{1-s}{2} \sum_{j=1}^{N-1} e^{- \kappa |j-i|} \iint_{\Gamma_j \times \Gamma_{j+1}} \frac{ ( w_j(x')^{-q} \eta_j(x') - w_j(y')^{-q} \eta_j(y'))^2 }{|x-y|^{n+s}} \, d \cH^{n-1}_x \, d \cH^{n-1}_y + C(1+M) \\
&\leq \frac{c_{n,s}}{100} \sum_{j=1}^{N-1} e^{- \kappa |j-i|} \int_{B_{r}^{n-1}} w_j^{-2q-1-s} \eta_j^2 + C_M(1+M).
\end{align*}
In particular, the coefficient $\frac{c_{n,s}}{100}$ is smaller than $\frac{1}{2}$ times the coefficient of the above left-hand side.
Hence, whenever $q \leq \frac{41}{20}$, the above inequality inequality implies the bound
\[
\frac{c_{n,s}}{100} \sum_{j=1}^{N-1} e^{- \kappa |j-i|} \int_{B_r^{n-1}} w_j^{-2q-1-s} \eta_j^2 \leq C_M(1+M). 
\]
This proves the desired property~\eqref{eqn:iterated-bound}.
To summarize, we have established the improvement ~\eqref{eqn:inductive-hypothesis} $\Rightarrow$ \eqref{eqn:iterated-bound}, showing that an $L^{-2q}$ estimate in a larger ball implies an $L^{-2q-1-s}$ estimate in a smaller ball, for $s \in ( \frac{1}{2},1)$ sufficiently close to $1$ and $q \leq \frac{41}{20}$.
Moreover, for $s \in ( \frac{9}{10}, 1)$, we have that
\[
2 \cdot 1 + 1+s > 2 \cdot \tfrac{3}{2}, \qquad 2 \cdot \tfrac{3}{2} + 1 + s > 2 \cdot \tfrac{41}{20}, \qquad 2\cdot \tfrac{41}{20} + 1+s > 6
\]
whereby applying this chain of improvements with $q \in \{ 1, \frac{3}{2}, \frac{41}{20} \}$ implies the desired $L^{-6}$ estimate in a smaller ball.
This completes the proof of our assertion.
\end{proof}

We now combine the $L^{-6}$ separation estimate of Proposition~\ref{prop:improve-the-exponent} with the results of~\cite{nonlocal-approximation}*{\S 3.5 and \S 4} to establish the structure of stable $s$-minimal sets described in Theorem~\ref{thm:sheet-happens}.
\begin{proof}[Proof of Theorem~\ref{thm:sheet-happens}]
The optimal sheet separation estimate is obtained as in~\cite{nonlocal-approximation}*{Propositions 3.5 and 3.6}, using the $L^{-6}$ estimate of Proposition~\ref{prop:improve-the-exponent} in place of their $L^{-2}$ estimate.
We retain the notation $w_i := (1-s)^{- \frac{1}{1+s}} (g_{i+1} - g_i)$ from the proof of Proposition~\ref{prop:improve-the-exponent}, so $w_i$ is uniformly equivalent to $\frac{g_{i+1} - g_i}{\sqrt{1-s}}$ as $s \uparrow 1$.
We choose a fixed $\theta \in (0, \frac{1}{64})$ so that the $L^{-6}$ estimate holds on $B^{n-1}_{8 \theta}$ in the fixed interior cylinders used below, after decreasing $\delta_*$ and increasing $s_*$ if necessary.

We first consider a pair with $g_i(0) = 0$ and $0 < g_{i+1}(0) \leq \frac{1}{8}$ and choose a fixed $\theta \in (0, \frac{1}{64})$ such that the $L^{-6}$ estimate holds on $B^{n-1}_{8 \theta}$.
We will also consider the neighboring sheets immediately outside the selected pair, namely $g_i - g_{i-1}$ and $g_{i+2} - g_{i+1}$, whenever they exist.

Now, on the one hand, if $g_i(0) - g_{i-1}(0) < \frac{1}{8}$, then $- \frac{1}{8} < g_{i-1}(0) < g_i(0) = 0$; likewise, if $g_{i+2}(0) - g_{i+1}(0) < \frac{1}{8}$, then $0 < g_{i+1}(0) < g_{i+2}(0) < \frac{1}{4}$.
Thus, the neighboring pair also satisfies the height condition required to apply the $L^{-6}$ estimate of Proposition~\ref{prop:improve-the-exponent} on its prescribed smaller $(n-1)$-dimensional ball.

On the other hand, if one pair of neighboring sheets has separation at least $\frac{1}{8}$ at $0$, then the assumed slope bound gives
\[
|D(g_{j+1} - g_j)| \leq |Dg_{j+1}| + |Dg_j| \leq 2 \delta. 
\]
Consequently, for $x' \in B^{n-1}_{1/2}$ and $\delta_* < \frac{1}{16}$, we obtain
\[
(g_{j+1} - g_j)(x') \geq (g_{j+1} - g_j)(0) - 2 \delta \, |x'| \geq \tfrac{1}{8} - \delta \geq \tfrac{1}{16}.
\]
In this case, the reciprocal of $w_j$ is uniformly bounded, with
\[
w_j^{-s} = (1-s)^{\frac{s}{1+s}} (g_{i+1} - g_i)^{-s} \leq 16^s (1-s)^{\frac{s}{1+s}} \leq 16.
\]
Thus, we can follow the computations of~\cite{nonlocal-approximation}*{(3.31) - (3.36)}, which produce a difference equation relating adjacent sheets.
We therefore obtain, in all cases,
\begin{equation}\label{eqn:estimates-on-the-wi}
\begin{split}
    |(1-s) \cL_{K_i} w_i| &\leq C( 1 + w_{i-1}^{-s} + w_i^{-s} + w_{i+1}^{-s}), \\
    \| (1-s) \cL_{K_i} w_i \|_{L^6(B^{n-1}_{4 \theta})} &\leq C, \\
    \sqrt{1-s} \, \| (g_{j+1} - g_j)^{-1} \|_{L^6(B^{n-1}_{4 \theta})} &\leq C, \qquad \text{for } \; j \in \{ i-1,i,i+1\}.
\end{split}
\end{equation}
The missing terms are omitted for $i=1$ and for $i=N-1$.
In particular, the right-hand side of the normalized difference equation for $g_{i+1} - g_i$ in~\cite{nonlocal-approximation}*{(3.31) - (3.35)}, including its truncation terms, has $L^6$ norm at most $C \sqrt{1-s}$.
We may therefore apply the inhomogeneous Harnack inequality for nonlocal operators~\cites{nonlocal-harnack , chaker-kim-weidner } as in the step preceding (3.36) in~\cite{nonlocal-approximation}, with $q = 6 > \frac{n-1}{1+s}$.
We therefore obtain, uniformly for $n \leq 8$ and $s$ sufficiently close to $1$,
\begin{equation}\label{eqn:little-harnack}
    \sup_{B_{\theta}^{n-1}} w_i \leq C (w_i(0) + 1),
\end{equation}
after dividing by $(1-s)^{\frac{1}{1+s}} \asymp \sqrt{1-s}$.

In the next step, we fix a point $x'_0 \in B^{n-1}_{\theta}$ and set $d = w_i(x'_0) \asymp \frac{(g_{i+1} - g_i)(x_0)}{\sqrt{1-s}}$.
It suffices to consider $d< \frac{\theta^2}{8}$.
For $d \leq r \leq \theta$, we translate $x_0$ to $(x'_0, g_i(x'_0))$ and dilate the configuration by $r^{-1}$: the slopes are unchanged, the Hessians are multiplied by $r$, and the heights of the selected rescaled sheets $\tilde{g}_i$ and $\tilde{g}_{i+1}$ over the origin are 
\[
0 \qquad \text{and} \qquad (1-s)^{\frac{1}{1+s}} r^{-1} d \leq (1-s)^{\frac{1}{1+s}} < \tfrac{1}{8}.
\]
Consequently, we can use the preceding argument in combination with the $L^{-6}$ estimate of Proposition~\ref{prop:improve-the-exponent} to each dilated configuration, in order to obtain
\begin{equation}\label{eqn:rescaled-wi-estimate}
    \sup_{B_{\theta r}^{n-1}(x'_0)} w_i \leq C(d+r), \qquad \text{for } \; d \leq r \leq \theta.
\end{equation}
We may therefore apply this discussion with $r = 2 \theta^{-1} |x'-x'_0|$ to obtain
\begin{equation}\label{eqn:wi-x-prime}
    w_i(x') \leq C \, |x'-x'_0|, \qquad d \leq |x' - x'_0| \leq \tfrac{1}{2} \theta^2.
\end{equation}
First, for $n \geq 7$, the $L^{-6}$ estimate of Proposition~\ref{prop:improve-the-exponent}, applied in the original cylinder, implies that
\begin{equation}\label{eqn:jacobian-computation}
C \geq \int_{B_{8 \theta}^{n-1}} w_i^{-6} \geq c \int_d^{\frac{1}{2} \theta^2} r^{n-8} \, dr \geq c \int_d^{\frac{1}{2} \theta^2} \frac{dr}{r} = c \log \frac{\theta^2}{2d} .
\end{equation}
Here, we computed as in~\cite{nonlocal-approximation}*{(3.37)}, and used $0 < r \leq \frac{1}{2} \theta^2 < 1$.
Thus, $d \geq \frac{1}{2} \theta^2 \exp( - \frac{C}{c}) > 0$, hence $\frac{(g_{i+1} - g_i)(x_0)}{\sqrt{1-s}} \asymp w_i(x'_0) = d \geq c' > 0$ proves the claimed separation estimate.

For $n=8$, the separation estimate follows from the conical modification of the above argument described in~\cite{nonlocal-approximation}*{Proposition 3.6}, which gains the one dimension needed in the bound~\eqref{eqn:jacobian-computation}.
This completes the proof of the first assertion.

Finally, given the separation estimate $\inf_{B_{r_*}^{n-1}} (g_{i+1} - g_i) \geq c \sqrt{1-s}$ for adjacent sheets, we can repeat the proofs of~\cite{nonlocal-approximation}*{Corollary 3.7 and Proposition 3.8}, using the nonlocal Schauder estimates of~\cite{fall}, to obtain the claimed uniform $C^{2,1/8}$ estimate,
\[
\| g_i \|_{C^{2,1/8}(B^{n-1}_{r_*/4})} \leq C.
\]
The same arguments bound the normalized fractional mean curvature of the selected sheet by $C \, \sqrt{1-s}$, and the comparison with the classical mean curvature as in Step 1 of the proof of~\cite{nonlocal-approximation}*{Lemma 3.9 and (3.32), (3.41)} gives the second estimate because $(1-s) \, |\log(1-s)| = o(\sqrt{1-s})$.
All steps above were carried out independently of the selected sheet $i$, the number $N$ of sheets, and the location of the vertex $x_{\circ}$, so the resulting bounds are independent of this data.
\end{proof}

\section{Proof of the main theorem}

We now explain how Theorem~\ref{thm:sheet-happens} can be combined with the results of~\cite{nonlocal-approximation}*{\S 3.5 and \S 4} and with the stable Bernstein theorem for minimal surfaces to prove Theorem~\ref{thm:stable-cones-s-1}.
Following Theorem~\ref{thm:sheet-happens}, we can argue as in~\cite{nonlocal-approximation}*{\S 3.7} to pass stability to the limiting lamination extracted from a sequence of stable $s_k$-minimal sets with $s_k \uparrow 1$.
Here, additional care is needed in carrying out the compactness and stability steps while the number of sheets may initially be unbounded.
\begin{lemma}\label{lemma:stability-lamination}
Consider an open set $U \subset \bR^n$ and a sequence of sets $E_k$ that are stationary and stable for the $s_k$-perimeter on every relatively compact subset of $U$, where $s_k \uparrow 1$.
Suppose that: 
\begin{enumerate}[(i)]
    \item The boundaries $\partial E_k$ are smooth and converge locally in the graphical $C^{2,\gamma}$ topology, for some $\gamma>0$, to a minimal lamination of $U$, with uniform local graphical bounds.
    \item In every fixed interior chart, adjacent approximating sheets are separated by at least $c \sqrt{1 -s_k}$, where $c>0$ may depend on the chart, but not on $k$.
\end{enumerate}
Then, the orientation cover of every leaf of the limiting minimal lamination is stable.
\end{lemma}
\begin{proof}
We denote the given lamination by $\cL$.  
We first observe that $\cL$ has locally Lipschitz charts: the differences of ordered minimal graphs satisfy uniformly elliptic linear equations, and the Harnack inequality combined with interior gradient estimates controls the $C^1$ difference of the graphs by their separation at a fixed transverse point.
Consequently, we can proceed as in Assertion 1 in the proof of~\cite{meeks-perez-ros}*{Theorem 1} to obtain locally bi-Lipschitz coordinate charts for $\cL$, hence $\cL$ has the Lipschitz regularity required to apply~\cite{meeks-perez-ros}*{Theorem 1 and Remark 1}.
These results imply the stability of limit leaves of $\cL$.

It remains to examine the isolated leaves of $\cL$.
Let $L$ be such a leaf, and we claim that
\begin{equation}\label{eqn:stability-inequality-L}
 \int_L |A_L|^2 \phi^2 \leq \int_L |\nabla^L \phi|^2, \qquad \text{for every } \; \phi \in C_c^{\infty}(L).   
\end{equation}
We notice that if $L$ is not two-sided, then letting $\pi: \tilde{L} \to L$ be its orientation double cover and taking any $\psi \in C^{\infty}_c(\tilde{L})$, we can apply~\eqref{eqn:stability-inequality-L}, by approximation, to the function
\[
\phi(x) = \Bigl( \sum_{\pi(y) = x} \psi(y)^2 \Bigr)^{\frac{1}{2}}.
\]
By Kato's inequality, we have $|\nabla \phi|^2 \leq \sum_{\pi(y) = x} |\nabla \psi(y)|^2$ almost everywhere.
Hence, combining this bound with~\eqref{eqn:stability-inequality-L} will imply the corresponding stability inequality for every $\psi \in C_c^{\infty}(\tilde{L})$.

We observe that the inequality~\eqref{eqn:stability-inequality-L}, for a given test function $\phi$, follows from the arguments of Proposition 3.1 and Lemmas 3.10-3.11 of~\cite{nonlocal-approximation}, as well as Step 3 in the proof of their Theorem 1.1.
Indeed, we can choose connected, relatively compact smooth domain $D \subset L$ containing $\textup{spt} \, \phi$ in its interior, so that a fixed normal tube over a neighborhood of $\bar{D}$ contains no other limiting local sheets.
For large $k$, normal projection of the $\partial E_k$ is a local diffeomorphism, by graphical convergence, which is proper by the convergence of tangent spaces.
Hence, this map is a local diffeomorphism of degree $N_k$.
On each fixed normal fiber, the hypothesis $(ii)$ gives $N_k \leq \frac{C}{\sqrt{1-s_k}}$.

Next, we choose a smooth ambient function $\eta$, compactly supported in this tube, which equals the pullback of $\phi$ near the approximating sheets.
This construction is well-defined even when $L$ is one-sided, using an even cutoff along the normal bundle.
To conclude, we apply the localized stability inequality of~\cite{nonlocal-approximation}*{Proposition 3.1}, multiply by $(1-s_k)$, and divide by $N_k$ as in case $(i)$ of Step 3 in the proof of Theorem 1.1 in~\cite{nonlocal-approximation}*{\S 3.7}.
The argument of Lemma 3.10 therein shows that the terms coming from the same local sheet converge to the two sides of the classical stability inequality~\eqref{eqn:stability-inequality-L}, multiplied by the same dimensional factor $c_n$ and uniform errors.
Moreover, all other potential terms are non-negative and can be discarded.

By~\cite{nonlocal-approximation}*{Lemma 3.11}, each ordered pair of distinct local sheets contributes at most $C(\phi) (1-s_k) (1 + |\log (1 - s_k)|)$.
Also, there are at most $N^2_k$ such pairs in this chart, with total contribution after averaging bounded by
\[
C(\phi) (1-s_k) N_k (1 + |\log(1-s_k)|) \leq C(\phi) \sqrt{1-s_k} (1 + |\log(1-s_k)|) \to 0.
\]
The remaining terms are bounded as in~\cite{nonlocal-approximation}*{Lemma 3.11}, proving~\eqref{eqn:stability-inequality-L} and the Lemma.
\end{proof}

Using Theorem~\ref{thm:sheet-happens} and Lemma~\ref{lemma:stability-lamination}, the above separation and regularity estimates imply a uniform curvature bound on the links, by the arguments of~\cite{nonlocal-approximation}*{\S 3}.
We obtain the following result.
\begin{corollary}\label{cor:curvature}
For $3\leq n\leq 8$, there are $s_* \in (0,1)$ and $C\in(0,\infty)$ such that for $s \in (s_*, 1)$, every stable $s$-minimal cone $E \subset \bR^n$ with vertex at the origin and $C^2$ link satisfies
\[
 \sup_{\partial E\cap\bS^{n-1}}|A_{\partial E}|\leq C,\qquad \text{for} \quad s\in(s_*,1).
\]
\end{corollary}
\begin{proof}
This result follows from applying White's blowup argument as in~\cite{nonlocal-approximation}*{\S 3.6 - 3.7}.
We argue by contradiction: if the assertion failed, there would exist $s_k \uparrow 1$, stable $s_k$-minimal cones $E_k$, and points $p_k \in \partial E_k \cap \bS^n$ where the curvature is maximal and satisfies $\lambda_k := |A_{\partial E_k}(p_k)| \to \infty$.
We dilate about $p_k$ by $\lambda_k$, so the vertices of the rescaled cones have distance $\lambda_k \to \infty$ from the origin and hence Theorem~\ref{thm:sheet-happens} applies; here, we use the conical structure of $E$ for $n=8$.

By the uniform local graphical estimates and the compactness argument of~\cite{nonlocal-approximation}*{\S 3.7}, these curvature blow-ups converge, after passing to a subsequence, to a minimal lamination $\cL$ of $\mathbb{R}^n$.
Moreover,~\cite{meeks-perez-ros}*{Theorem 1 and Remark 1}, applied as in~\cite{nonlocal-approximation}*{Lemmas 3.10 and 3.11}, shows that the leaf $L$ through the origin is complete, satisfies $|A_L|(0) = 1$, and has stable orientation cover, by Lemma~\ref{lemma:stability-lamination}.
For $n \leq 7$, this would contradict the stable Bernstein theorem.

For $n=8$, after passing to a subsequence, we can arrange that $p_k \to e$, for a constant vector $e$, so the radial tangent directions of the original cones converge to $e$.
Hence, $e$ is tangent to the limit leaf $L$ through the origin everywhere, and since $L$ is complete, we can write $L = L_0 \times \bR e$.
The stability of $L$ implies that $L_0 \subset e^{\perp} \cong \bR^7$ is also stable, contradicting the stable Bernstein theorem in $\mathbb{R}^7$ due to $|A_L(0)| = 1$.
Hence, $\sup_{\partial E \cap \bS^{n-1}} |A_{\partial E}| \leq C$ and our assertion is proved.
\end{proof}

Next, we extend the graphical decomposition of~\cite{nonlocal-approximation}*{Lemma 4.5} to hypersurfaces of $\mathbb{S}^{n-1}$, for $3 \leq n\leq 8$, using the lamination theory of~\cite{meeks-perez-ros} and the results of~\cites{catino-mastrolia-roncoroni , wang-bernstein }.
\begin{lemma}\label{lemma:limiting-lamination}
For $3 \leq n \leq 8$ and $\alpha \in (0,1)$, let $S_k \subset \bS^{n-1}$ be a sequence of oriented, embedded, closed hypersurfaces satisfying uniform local $C^{2,\alpha}$ graphical bounds, with mean curvatures
\[
\| H[S_k] \|_{L^{\infty}(S_k)} \to 0 \qquad \text{as} \qquad k \to \infty.
\]
Then, after passing to a subsequence, there is a connected, closed, embedded minimal hypersurface $\Sigma \subset \bS^{n-1}$, with unit normal $\nu$, such that
\[
S_k = \bigcup_{i=1}^{N_k} \Biggl\{ \frac{x + u_{k,i}(x) \, \nu(x)}{\sqrt{1 + u_{k,i}(x)^2}} : x \in \Sigma \Biggr\}
\]
where $u_{k,i} \in C^2(\Sigma)$ for all $i \in \{ 1, \dots,N_k \}$ and
\[
u_{k,1} < \dots < u_{k,N_k}, \qquad \max_i \| u_{k,i} \|_{C^{2,\beta}(\Sigma)} \to 0, \quad \text{for every} \quad 0<\beta<\alpha.
\]
In particular, we do not assume any area bound or any uniform bound on $N_k$.
\end{lemma}
\begin{proof}
We follow the compactness argument in Step 1 of the proof of~\cite{nonlocal-approximation}*{Lemma 4.5}, which is independent of the dimension under the assumed uniform $C^{2,\alpha}$ bounds. 
The uniform graphical bounds and Arzel\`a-Ascoli give, after passing to a subsequence, local $C^{2,\beta}$ convergence for every $\beta<\alpha$.
Because the mean curvatures of the $S_k$ tend to zero, the limiting local sheets are minimal.
Embeddedness and the strong maximum principle imply that the leaves are either disjoint or agree, so their union defines a minimal lamination $\cL$ of $\bS^{n-1}$.
The local regularity and completeness arguments used in the proof of Lemma~\ref{lemma:stability-lamination} apply in $\bS^{n-1}$, whereby $\cL$ has complete leaves and satisfies the Lipschitz regularity required to apply~\cite{meeks-perez-ros}*{Theorem 1 and Remark 1}.

Thus, the orientation cover of every limit leaf is stable.
For $3 \leq n \leq 8$, such a cover would be a complete, orientable, stable minimal immersion of dimension $n-2 \leq 6$ into $\bS^{n-1}$, contradicting~\cite{catino-mastrolia-roncoroni}*{Corollary 1.3} and~\cite{wang-bernstein}*{Corollary 1.6}.

We can modify Step 2 in the proof of~\cite{nonlocal-approximation}*{Lemma 4.5} accordingly: if a limiting image were not a compact embedded hypersurface, we could repeat the construction of a positive Jacobi field therein to obtain orientable stable immersed minimal intrinsic balls of arbitrarily large radius, with the same uniform local graphical bounds.

Applying the pointed compactness argument of Step 1 once more, and passing the unit normals to the limit, we would obtain an orientable minimal immersion as the radii tend to infinity.
Because the intrinsic distances from the base points to the boundaries tend to infinity, such a limiting immersion is complete.
Moreover, since every compactly supported test function is contained in a smoothly converging interior region, the immersion is stable.
For $3 \leq n\leq 8$, this gives the same contradiction to~\cite{catino-mastrolia-roncoroni}*{Corollary 1.3} and~\cite{wang-bernstein}*{Corollary 1.6} as discussed above.

We conclude that $\cL$ has no limit leaves, so it consists of finitely many compact embedded leaves.
Indeed, a non-compact leaf, or infinitely many distinct leaves, would produce an accumulating sequence of local sheets, and hence a limit leaf.
Moreover, Frankel's theorem~\cite{frankel} shows that any two distinct compact leaves would intersect, producing a contradiction.
Therefore, $\cL$ consists of a single connected, closed, embedded minimal hypersurface $\Sigma$, which is two-sided.
%For $n=3$, the complete leaves are great circles in $\bS^2$, whose intersection gives the same conclusion directly.

We conclude as in Step 3 of the proof of~\cite{nonlocal-approximation}*{Lemma 4.5}: because $\Sigma$ is two-sided, for large $k$, the hypersurface $S_k$ lies in a fixed normal tube around $\Sigma$, and normal projection is a local diffeomorphism, hence a finite covering due to compactness.
By the ordering of the points in each normal fiber, this finite cover is trivial, hence each $S_k$ consists of globally ordered normal graphs over $\Sigma$, which converge uniformly to $0$ in every $C^{2,\beta}$ norm.
This proves our assertion.
\end{proof}

The lamination argument employed above, using~\cite{wang-bernstein}*{Corollary 1.6}, applies verbatim to the study of stable classical minimal cones in $\mathbb{R}^8$.
Consequently, it allows us to remove the density assumption in Edelen's compactness result~\cite{edelen}*{Theorem 5.1}.
We obtain the following.
\begin{theorem}\label{thm:density-classical-cones}
Let $\cC$ denote the collection of embedded, stable minimal hypercones in $\mathbb{R}^8$, with smooth links.
There exists~$C\in(0,\infty)$ such that every $\mathbf{C} \in \cC$ satisfies
\[
\sup_{\Sigma} |A_{\Sigma}| \leq C \quad \text{and} \quad \cH^6(\Sigma) \leq C, \qquad \text{where} \quad \Sigma = \mathbf{C} \cap \bS^7.
\]
Moreover, every sequence $\{ \mathbf{C}_k \} \subseteq \cC$ admits a subsequence whose links converge smoothly, with multiplicity one, to the link of another cone in $\cC$.
\end{theorem}

In the next step, we can repeat~\cite{nonlocal-approximation}*{Lemma 4.8} independently of the ambient dimension to obtain a uniform bound on the number of sheets, using the separation estimate of Theorem~\ref{thm:sheet-happens}.
\begin{lemma}\label{lemma:bound-the-sheets}
Let $\Sigma^{n-2} \subset \bS^{n-1}$ be a connected, closed, embedded minimal hypersurface with unit normal $\nu$.
Suppose that the conical normal graphs over $\Sigma$
\[
\Gamma_{k,i} := \{ rx + r \, u_{k,i}(x) \nu(x) : r > 0, \; x \in \Sigma \}, \qquad 1\leq i \leq N_k,
\]
are ordered by $u_{k,1} < \cdots < u_{k,N_k}$ and satisfy $\max_i \| u_{k,i} \|_{C^{2,\gamma}(\Sigma)} \to 0$ and
\[
\sup_i \sup_{\Gamma_{k,i} \cap \bS^{n-1}} |H[\Gamma_{k,i}]| \leq C \sqrt{1-s_k} \qquad \text{and} \qquad u_{k,i+1} - u_{k,i} \geq c \sqrt{1-s_k} \, .
\]
Then, it holds that
\[
\| u_{k,N_k} - u_{k,1} \|_{L^{\infty}(\Sigma)} \leq C' \sqrt{1-s_k} \qquad \text{and} \qquad N_k \leq C'
\]
for a constant $C'$ that is independent of $N_k$ and $k$.
\end{lemma}
\begin{proof}
Let $\cH(u)$ denote the classical mean curvature operator for the conical normal graph, evaluated at $x + u(x) \nu(x)$, with sign chosen so that $D \cH(0) = \Delta_{\Sigma} + |A_{\Sigma}|^2 + n-2$ is the Jacobi operator.
Then, taking the mean-value linearization of $\cH$ along the segment joining two graphical sheets $u_{k,i}$ and $u_{k,j}$, we obtain uniformly elliptic operators $L^{(k)}_{i,j} = \int_0^1 D \cH ( u_{k,j} + t(u_{k,i} - u_{k,j})) \, dt$ with
\begin{equation}\label{eqn:linearized-operators}
L^{(k)}_{i,j} (u_{k,i} - u_{k,j}) = \cH( u_{k,i}) - \cH (u_{k,j}), \qquad L^{(k)}_{i,j} \to \Delta_{\Sigma} + |A_{\Sigma}|^2 + (n-2).
\end{equation}
The convergence of the operators $L^{(k)}_{i,j}$ is understood at the level of coefficients.

We now specialize this computation to the two extreme sheets, taking $(i,j) = (N_k,1)$ and considering the difference $u_{k,N_k} - u_{k,1}$.
If $\frac{1}{\sqrt{1-s_k}}\| u_{k,N_k} - u_{k,1} \|_{L^{\infty}} \to \infty$, we normalize $h_k := \frac{u_{k,N_k} - u_{k,1}}{\| u_{k,N_k} - u_{k,1} \|_{L^{\infty}}}$.
Then, applying standard elliptic estimates as in~\cite{gilbarg-trudinger}, together with the Arzel\`a-Ascoli theorem, produces a non-negative limit $h$ satisfying
\[
\| h\|_{L^{\infty}} = 1, \qquad ( \Delta_{\Sigma} + |A_{\Sigma}|^2 + (n-2)) h = 0 .
\]
Because $\int_{\Sigma} \Delta_{\Sigma}h =0$, we obtain a contradiction after integrating over $\Sigma$.
Thus, there is a constant $C'$, independent of $k$ and $N_k$, for which
\[
\| u_{k,N_k} - u_{k,1} \|_{L^{\infty}(\Sigma)} \leq C' \sqrt{1-s_k}.
\]
On the other hand, $u_{k,i+1} - u_{k,i} \geq c \sqrt{1-s_k}$ for each $i$.
Combining these properties, we obtain
\[
(N_k - 1) c \leq (1-s_k)^{- \frac{1}{2}} \| u_{k,N_k} - u_{k,1} \|_{L^{\infty}(\Sigma)} \leq C'
\]
which implies that $N_k \leq C''$ is uniformly bounded.
This completes the proof of our assertion. 
\end{proof}

Finally, we obtain the higher-dimensional analogue of~\cite{nonlocal-approximation}*{Proposition 4.6}.
\begin{corollary}\label{corollary:my-favorite-russian}
Let $3 \leq n \leq 8$ and suppose that $s_k \uparrow 1$ and that $E_k \subset \mathbb{R}^n$ are non-flat stable $s_k$-minimal cones with vertex at the origin and $C^2$ links.
After passing to a subsequence, there exist a connected, closed, embedded minimal hypersurface $\Sigma \subset \mathbb{S}^{n-1}$, with unit normal $\nu$, an integer $N \geq 1$, and globally ordered functions $u_{k,1} < \cdots < u_{k,N}$ on $\Sigma$, such that
\[
\partial E_k \setminus \{ 0 \} = \bigcup_{i=1}^N \Gamma_{k,i}, \qquad \Gamma_{k,i} = \{ rx + r \, u_{k,i}(x) \nu(x) : r >0, \; x \in \Sigma\}. 
\]
The cone $C(\Sigma)$ is stable away from the origin, and, for every given $\gamma \in (0, \frac{1}{16})$, we have that $$\max_i \| u_{k,i} \|_{C^{2,\gamma}(\Sigma)} \to 0.$$
Moreover, there are constants $0<c<C<\infty$, independent of $k$ and $i$, such that
\begin{align}
    \| u_{k,N} - u_{k,1} \|_{L^{\infty}(\Sigma)} &\leq C \sqrt{1-s_k}, \qquad u_{k,i+1} - u_{k,i} \geq c \sqrt{1-s_k} \qquad \text{on } \; \Sigma, \label{eqn:limiting-separation}
    \intertext{and}
    \max_i \sup_{\Gamma_{k,i} \cap \bS^{n-1}} |H[\Gamma_{k,i}]| &\leq C \sqrt{1-s_k}. \label{eqn:mean-curvature-bound}
\end{align}
In particular, the boundaries converge graphically to $C(\Sigma)$ in $C^{2,\gamma}$ on every compact annulus, with fixed multiplicity $N$.
The width and separation estimates extend by homogeneity.
\end{corollary}
\begin{proof}
By Theorem~\ref{thm:sheet-happens} and Corollary~\ref{cor:curvature}, the links $\partial E_k \cap \bS^{n-1}$ have uniform local $C^{2,1/8}$ graphical bounds and classical mean curvature bounded by $C \sqrt{1-s_k}$.
The spherical compactness result of Lemma~\ref{lemma:limiting-lamination} therefore gives, after passing to a subsequence, globally ordered normal graphs over a connected, closed, embedded minimal hypersurface $\Sigma \subset \bS^{n-1}$ with unit normal $\nu$, converging uniformly to zero in $C^{2,\gamma}$ for any fixed $\gamma \in (0, \frac{1}{16})$.
Since $E_k$ is conical, we may therefore write
\begin{align*}
    \partial E_k \cap \bS^{n-1} &= \bigcup_i \left\{ \frac{x + u_{k,i}(x) \nu(x)}{\sqrt{1 + u_{k,i}(x)^2}} : x \in \Sigma \right\}, \qquad \max_i \| u_{k,i} \|_{C^{2,\gamma}(\Sigma)} \to 0, \\
    \partial E_k \setminus \{ 0 \} &= \bigcup_i \{ rx + r \, u_{k,i}(x) \nu(x) : r>0, \; x \in \Sigma \}.
\end{align*}
By homogeneity, the cones $\partial E_k$ graphically on compact annuli to $C(\Sigma)$.
Moreover, the property $\max_i \| u_{k,i} \|_{C^{2,\gamma}(\Sigma)} \to 0$ implies that, on every fixed radial interval $0 < a \leq r \leq b < \infty$, we have
\[
\max_i \| (r,x) \mapsto r \, u_{k,i}(x) \|_{C^{2,\gamma}( [a,b] \times \Sigma)} \leq C(a,b) \, \max_i \| u_{k,i} \|_{C^{2,\gamma}(\Sigma)} \to 0.
\]
Consequently, the conical sheets $\partial E_k$ converge to $C(\Sigma)$ graphically in $C^{2,\gamma}$ on every compact subset of $U = \bR^n \setminus \{ 0 \}$, with uniform local graphical bounds.

We observe that the hypotheses $(i), (ii)$ of Lemma~\ref{lemma:stability-lamination} are satisfied for the one-leaf minimal lamination of $U = \bR^n \setminus \{ 0 \}$ defined by the smooth embedded minimal hypersurface $C(\Sigma) \setminus \{ 0 \}$.
The first property is the conclusion of the above discussion, while the second property follows from Theorem~\ref{thm:sheet-happens} after a local rescaling in a small Euclidean chart so that, after translating, rotating and dilating this chart by $r^{-1}$, the transformed sets remain stationary and stable, while the limiting cone is a graph with slope and Hessian smaller than $\frac{1}{2} \delta_*$.
Thus, applying Theorem~\ref{thm:sheet-happens} and rescaling back to the original coordinates shows that adjacent sheets are separated by at least $c \, r \sqrt{1 - s_k}$, where the constant $c \,r$ depends only on the chart, but not on $k$. 
Each compact region of the limiting cone $C(\Sigma)$ is covered by a finite collection of such charts, hence $(ii)$ holds locally throughout $U$.

To summarize, we can apply Lemma~\ref{lemma:stability-lamination} on $U = \bR^n \setminus \{ 0 \}$ to deduce that $C(\Sigma)$ is stable away from the origin in $\bR^n$.
Moreover, Theorem~\ref{thm:sheet-happens} gives uniform sheet separation estimates and the uniform bound on the mean curvatures of the homogeneous normal graphs of the $u_{k,i}$.
These properties allow us to apply Lemma~\ref{lemma:bound-the-sheets} to obtain~\eqref{eqn:limiting-separation} as well as a uniform bound $N_k \leq C'$, and passing to a further subsequence gives $N_k = N$.
This proves our assertion.
\end{proof}

We recall the following standard result, which follows by arguments similar to those of~\cite{nonlocal-approximation}*{\S 4.3 - 4.5}, described also in~\cite{nonlocal-approximation}*{\S 1.4}.
\begin{lemma}\label{lemma:limiting-problem}
In the setting of Corollary~\ref{corollary:my-favorite-russian}, after passing to a subsequence, the functions $\frac{u_{k,i+1} - u_{k,i}}{\sqrt{1-s}}$ converge to positive functions $v_i$ on $\Sigma$, for $1 \leq i < N$, which satisfy
\begin{equation}\label{eqn:jacobi-vi}
    (\Delta_{\Sigma} + |A_{\Sigma}|^2 + n-2) v_i \leq  4 v_i^{-1}.
\end{equation}
If $N=2$, this relation holds with equality.
Moreover, for every smooth function $\phi$, it holds that
\begin{equation}\label{eqn:stability-inequality}
    \int_{\Sigma} |A_{\Sigma}|^2 \phi^2 + \frac{8}{N} \sum_{i=1}^{N-1} \int_{\Sigma} \frac{\phi^2}{v_i^2} \leq \int_{\Sigma} |\nabla^{\Sigma} \phi|^2 + \frac{(n-3)^2}{4} \int_{\Sigma} \phi^2.
\end{equation}
\end{lemma}
\begin{proof}
We recall the general geometric expansion of~\cite{nonlocal-approximation}*{Lemma 3.9}, which expresses the classical mean curvature $\cH$ of the conical normal graphs in terms of the distances to the other sheets, up to an error of order $O( (1-s)^{\frac{1+\beta}{2}})$.
We apply this result in fixed local charts on compact annular regions of $C(\Sigma)$, where the uniform graphical convergence obtained in Corollary~\ref{corollary:my-favorite-russian} gives
\[
\textup{dist} (x + u_{k,i}(x) \nu(x), \Gamma_{k,j}) = |u_{k,j}(x) - u_{k,i}(x)| \, (1 + o(1))
\]
uniformly in $x,i,j$.
Indeed, the relative error is $o(1)$, as seen by integrating the signed distance to $\Gamma_{k,j}$ along the normal fiber between $u_{k,i}(x)$ and $u_{k,j}(x)$.
Consequently, since $N$ is fixed and the separation is bounded below by $c \sqrt{1-s_k}$, the expansion of the mean curvature $\cH(u_{k,i})$ becomes
\[
\cH (u_{k,i}) = 2 (1-s_k) \sum_{j \neq i} \frac{(-1)^{i-j}}{u_{k,j} - u_{k,i}} + o ( \sqrt{1-s_k}).
\]
Subtracting the equations of adjacent sheets and applying the standard mean value linearization~\eqref{eqn:linearized-operators} with $(i,j) = (i+1,i)$, also employed in Lemma~\ref{lemma:bound-the-sheets}, produces the equation
\[
\cH(u_{k,i+1}) - \cH (u_{k,i}) = L^{(k)}_i ( u_{k,i+1}- u_{k,i}), \qquad L^{(k)}_i \to \Delta_{\Sigma} + |A_{\Sigma}|^2 + (n-2),
\]
converging to the Jacobi equation.
Let $v_{k,i} := \frac{u_{k,i+1} - u_{k,i}}{\sqrt{1-s}}$ denote the normalized separation, which satisfies $c \leq v_{k,i} \leq C$ and $\| L^{(k)}_i v_{k,i} \|_{L^{\infty}(\Sigma)} \leq C$ in view of Lemma~\ref{lemma:bound-the-sheets}.
Hence, standard elliptic estimates~\cite{gilbarg-trudinger} imply that $\| v_{k,i} \|_{W^{2,p}(\Sigma)} \leq C_p$ for every $1<p<\infty$.

The interaction between the selected pair of sheets contributes a term $4 v_{k,i}^{-1}$, and arguing as in~\cite{nonlocal-approximation}*{\S 4.3 - 4.5} shows that each exterior tail term is an alternating sum that begins with a negative term and decreasing subsequence contributions.
Combining these observations, we obtain
\begin{equation}\label{eqn:Lki-vki}
L^{(k)}_i v_{k,i} \leq 4 v_{k,i}^{-1} + o(1).
\end{equation}
Using the uniform bounds on the $v_{k,i}$, we deduce that a subsequence of these functions converges strongly in $C^{1,\beta}$, for every $\beta<1$, and weakly in $W^{2,p}$, to a function $v_i$ satisfying $c \leq v_i \leq C$.
Therefore, the convergence $L^{(k)}_i \to \Delta_{\Sigma} + |A_{\Sigma}|^2 + (n-2)$, as uniformly elliptic operators, can be combined with the inequality~\eqref{eqn:Lki-vki} to establish the property~\eqref{eqn:jacobi-vi}.
When $N=2$, there are no exterior tail terms and~\eqref{eqn:jacobi-vi} holds with equality.

To obtain the stability condition~\eqref{eqn:stability-inequality}, we first observe, as in Lemma~\ref{lemma:stability-lamination}, that Lemmas 3.10-3.11 of~\cite{nonlocal-approximation} apply to normal graphs over a non-flat hypersurface, where we test the localized stability formula of~\cite{nonlocal-approximation}*{\S 3.1} with the same scalar normal speed on every sheet, pulled back from a function $\eta$ on $C(\Sigma)$.
The contributions from the same sheet converge, with a common positive dimensional factor $c_n$, to the term
\[
N c_n \int_{C(\Sigma)} ( |\nabla^{\Sigma} \eta|^2 - |A_{C(\Sigma)}|^2 \eta^2).
\]
As in~\cite{nonlocal-approximation}*{\S 4.3 - 4.5}, the mixed terms tend to zero, while for each adjacent pair, rescaling by its normal separation gives
\[
(1 - s_k) \int_{\Gamma_{k,i+1}} \frac{d \cH^{n-1}_y}{|rx + r u_{k,i}(x) \nu(x) - y|^{n+s_k}} = (c_n + o(1)) \frac{1-s_k}{r^2 (u_{k,i+1} - u_{k,i})^2}.
\]
After rescaling adjacent pairs by the normal separation and passing to the limit, the two tangent planes coincide, while the outward normals are opposite and hence have squared normal difference tending to $4$.
Including both orderings of this pair therefore produces the coefficient $8 c_n$, and discarding all other non-negative potential terms as in~\cite{nonlocal-approximation}*{Proposition 4.11}.

Finally, we apply these considerations to a test function $\eta(r,x) = \zeta(r) \phi(x)$ on $C(\Sigma)$, where $\zeta$ is supported in a fixed compact radial interval.
Consequently, after sending $k \to \infty$, using $|A_{C(\Sigma)}|^2 = r^{-2} |A_{\Sigma}|^2$, and separating radial and angular integrals, we obtain
\begin{align*}
    \int_{\Sigma} |A_{\Sigma}|^2 \phi^2 + \frac{8}{N} \sum_i \int_{\Sigma} \frac{\phi^2}{v_i^2} \leq \int_{\Sigma} |\nabla^{\Sigma} \phi|^2 + \frac{\int_0^{\infty} r^{n-2} |\zeta'|^2 \, dr}{\int_0^{\infty} r^{n-4} \zeta^2 \, dr} \int_{\Sigma} \phi^2.
\end{align*}
Taking the infimum of the radial quotient on the right-hand side recovers the Hardy threshold $(\frac{n-3}{2})^2$ and proves the desired inequality~\eqref{eqn:stability-inequality}.
\end{proof}

\subsection{Completion of the proof}\label{sec:completion}

We now apply these observations to obtain Theorems~\ref{thm:stable-cones-s-1} and~\ref{thm:eight-dimensions}.

\begin{proof}[Proof of Theorems~\ref{thm:stable-cones-s-1} and~\ref{thm:eight-dimensions}]
Let $3 \leq n \leq 8$, and suppose that $s_k \uparrow 1$ and that $E_k \subset \bR^n$ are non-flat stable $s_k$-minimal cones with $C^2$ links.
By Corollary~\ref{corollary:my-favorite-russian}, after passing to a subsequence, the boundaries $\partial E_k$ converge graphically on compact annuli, with fixed multiplicity $N$, to a stable classical minimal cone $C(\Sigma)$, for $\Sigma \subset \bS^{n-1}$ a connected, closed, embedded minimal hypersurface.

Suppose first that $N \geq 2$.
By Lemma~\ref{lemma:limiting-problem}, the separations between adjacent sheets, normalized by $\sqrt{1-s_k}$, converge to positive functions $\{ v_i \}_{i=1}^{N-1}$ that satisfy~\eqref{eqn:jacobi-vi} and~\eqref{eqn:stability-inequality} for every smooth function $\phi$ on $\Sigma$.
Multiplying the inequality~\eqref{eqn:jacobi-vi} by $v_i^{-1}$ and integrating by parts gives
\begin{equation}\label{eqn:stability-contradiction-I}
    (n-2) \, |\Sigma| + \int_{\Sigma} |A_{\Sigma}|^2 + \int_{\Sigma} |\nabla^{\Sigma} \log v_i|^2 \leq 4 \int_{\Sigma} v_i^{-2}.
\end{equation}

Moreover, choosing $\phi=1$ in the stability inequality~\eqref{eqn:stability-inequality} yields
\begin{equation}\label{eqn:stability-contradiction-II}
    \int_{\Sigma} |A_{\Sigma}|^2 + \frac{8}{N} \sum_{i=1}^{N-1} \int_{\Sigma} v_i^{-2} \leq \frac{(n-3)^2}{4} |\Sigma|.
\end{equation}
We sum the inequality~\eqref{eqn:stability-contradiction-I}, multiply by $\frac{2}{N}$, and combine the result with~\eqref{eqn:stability-contradiction-II}.
After discarding non-negative terms, this produces the inequality $\frac{2(N-1)}{N}(n-2) \leq ( \frac{n-3}{2})^2$.
For $3 \leq n \leq 7$, the bounds $n-2> \frac{(n-3)^2}{4}$ and $\frac{2(N-1)}{N} \geq 1$ for $N \geq 2$ produce a contradiction, hence $N=1$.
For $n=8$, the same computation shows that $N \geq 3$ is impossible, hence $N \in \{ 1,2\}$.

For $3 \leq n \leq 7$, the stable minimal cone $C(\Sigma) \subset \bR^n$ must be a hyperplane, and $N=1$.
After a rotation, the graphical representation of~\ref{corollary:my-favorite-russian} therefore expresses each $\partial E_k$, for sufficiently large $k$, as a single globally defined one-homogeneous graph over $\mathbb{R}^{n-1}$, smooth away from the origin.
By~\cite{flatness-farina-valdinoci}, each $E_k$ is then a half-space, contradicting the assumption that the cones are non-flat.

We now consider $n=8$.
If $N=1$, the same argument as above applies to show $C(\Sigma)$ cannot be a hyperplane, hence it must be a non-flat stable classical minimal cone.
In this case, the claimed uniform compactness follows from the compactness Theorem~\ref{thm:density-classical-cones} for classical minimal cones.
For $N=2$, Lemma~\ref{lemma:limiting-problem} shows that~\eqref{eqn:jacobi-vi} holds with equality for the function $v_1$.

On the other hand, we obtain the inequality
\[
\tfrac{1}{2} \Delta_{\Sigma} |A_{\Sigma}|^2 = |\nabla A_{\Sigma}|^2 + 6 \, |A_{\Sigma}|^2 - |A_{\Sigma}|^4 \geq |\nabla |A_{\Sigma}||^2 + 6 \, |A_{\Sigma}|^2 - |A_{\Sigma}|^4
\]
by applying Simons' identity and Kato's inequality.
Integrating this expression, we arrive at
\begin{equation}\label{eqn:integrating-simons-kato}
6 \int_{\Sigma} |A_{\Sigma}|^2 + \int_{\Sigma} |\nabla |A_{\Sigma}||^2 \leq \int_{\Sigma} |A_{\Sigma}|^4.
\end{equation}
Next, we multiply~\eqref{eqn:jacobi-vi} by $v_1^{-1} |A_{\Sigma}|^2$, integrate by parts, and complete the square, to find
\begin{align*}
    6 \int_{\Sigma} |A_{\Sigma}|^2 + \int_{\Sigma} |A_{\Sigma}|^4 &= 4 \int_{\Sigma} v_1^{-2} |A_{\Sigma}|^2 + 2 \int_{\Sigma} |A_{\Sigma}| \, v_1^{-1} \la \nabla |A_{\Sigma}|, v_1 \rg - \int_{\Sigma} |A_{\Sigma}|^2 v_1^{-2} |\nabla v_1|^2 \\
    &\leq 4 \int_{\Sigma} v_1^{-2} |A_{\Sigma}|^2 + \int_{\Sigma} |\nabla |A_{\Sigma}||^2.
\end{align*}
Combining this inequality with~\eqref{eqn:integrating-simons-kato}, we arrive at
\begin{equation}\label{eqn:times-three-inequality}
12 \int_{\Sigma} |A_{\Sigma}|^2 \leq 4 \int_{\Sigma} v_1^{-2} |A_{\Sigma}|^2, \qquad \implies \qquad 3 \int_{\Sigma} |A_{\Sigma}|^2 \leq \int_{\Sigma} v_1^{-2} |A_{\Sigma}|^2.
\end{equation}
On the other hand, a standard smooth approximation allows us to test the inequality~\eqref{eqn:stability-inequality} with $\phi_{\ve} = (|A_{\Sigma}|^2+\ve^2)^{\frac{1}{2}}$, and then send $\ve \downarrow 0$.
Combined with~\eqref{eqn:integrating-simons-kato}, this implies that
\[
4 \int_{\Sigma} v_1^{-2} |A_{\Sigma}|^2 \leq \int_{\Sigma} |\nabla |A_{\Sigma}||^2 + \frac{25}{4} \int_{\Sigma} |A_{\Sigma}|^2 - \int_{\Sigma} |A_{\Sigma}|^4 \leq \frac{1}{4} \int_{\Sigma} |A_{\Sigma}|^2.
\]
In view of the inequality~\eqref{eqn:times-three-inequality}, the latter bound forces $|A_{\Sigma}| \equiv 0$.
Since $\Sigma$ is connected, closed, and embedded, we conclude that it is an equator, and the boundaries $\partial E_k$ converge with multiplicity $N=2$ to a hyperplane through the origin.
This completes the proof.
\end{proof}

\bibliography{ref}

\end{document}